\documentclass [11pt]{article}
\usepackage{graphicx}
\usepackage{sfmath}

\begin{document}

\title{Stability of homogeneous bundles on ${\bf P}^3$}
\author{Elena Rubei}
\date{}
\maketitle

\def\thefootnote{}
\footnotetext{ \hspace*{-0.36cm}
{\bf Address}: Dipartimento di Matematica ``U.Dini'', Viale Morgagni 
67/A, c.a.p.  50134 Firenze, Italia. 

{\bf E-mail address:} rubei@math.unifi.it }



\newtheorem{theorem}{Theorem}
\newtheorem{lemma}[theorem]{Lemma}
\newtheorem{prop}[theorem]{Proposition}
\newtheorem{rema}[theorem]{Remark}
\newtheorem{cor}[theorem]{Corollary}
\newtheorem{defin}[theorem]{Definition}
\newtheorem{notat}[theorem]{Notation}
\newtheorem{noterec}[theorem]{Notation and recalls}
\newtheorem{exam}[theorem]{Example}

\vspace*{-0.5cm}

\begin{abstract} We study the stability and the simplicity 
of some homogeneous bundles on ${\bf P}^3$ by using the quiver
 associated to homogeneous bundles introduced by Bondal and Kapranov 
in \cite{B-K}.
In particular we show that the homogeneous bundles on ${\bf P}^3$ 
whose quiver support is a parallelepiped 
or a classical staircase are stable. 
For instance the bundles $E$ whose minimal
free resolution is of the kind
 $$0 \rightarrow S^{\lambda_1, \lambda_2, \lambda_3 } V (t) \rightarrow 
S^{\lambda_1+s, \lambda_2, \lambda_3 } V (t+s) 
\rightarrow E \rightarrow 0 $$ are stable.
\end{abstract}

\section{Introduction}

We examine some homogeneous bundles on  ${\bf P}^3$
whose minimal free resolutions are of a particular kind and we study
their simplicity and stability. The main tool we use is  the quiver 
associated to homogeneous bundles  introduced by Bondal 
and Kapranov in \cite{B-K}.

Quivers allow us to handle well and ``to make explicit'' the homogeneous 
subbundles of a homogeneous bundle $E$
and, by  Rohmfeld's criterion (see \cite{Rohm}),  
$E$ is semistable if and only if the slope   of every subbundle associated 
to a subrepresentation of the quiver representation of $E$ is less or equal 
than the slope of $E$; so quivers and representations of quivers associated
to homogeneous bundles are particularly suitable for the study of stability.

In this paper we focus on homogeneous bundles on ${\bf P}^3$
such that the support of their 
quiver representation (which we will call ``quiver support'' for the 
sake of brevity) 
is as simple as possible. In particular  
we prove (Theorems \ref{parstable} and \ref{stairstable}) 
that the homogeneous bundles whose quiver  support 
is a {\bf parallelepiped} are stable (if they are not tensor product of 
an $SL(V)$-representation and  ${\cal O}(t)$ for some $t$) and also the
 bundles whose quiver support 
has the form of a {\bf classical staircase} (see the third figure in \S4) 
are stable. 

One can easily prove that, 
 if $E$ is a
homogeneous vector bundle on ${\bf P}^3={\bf P}(V)$  there exists 
a minimal free resolution of $E$  $$ 0 \rightarrow \oplus_q
{\cal O}(-q) \otimes_{\bf C} A_{q} \rightarrow \oplus_q {\cal
O}(-q) \otimes_{\bf C} B_{q}  \rightarrow \oplus_q {\cal
O}(-q) \otimes_{\bf C} C_{q}\rightarrow E \rightarrow 0 $$ 
with $A_q, B_q, C_q$ $SL(V)$-representations and maps $SL(V)$-invariant. 

In terms of their minimal free resolutions the results on stability of the 
homogeneous bundles such that the
 support of their quiver representations is a parallelepiped 
or a classical staircase,  can be 
restated respectively in the following theorems:

\begin{theorem} \label{parstablebis}
Let $E$ be a homogeneous bundle on ${\bf P}^3={\bf P}(V) $ whose 
minimal free resolution is one of the following 
for some $\lambda_1, \lambda_2, \lambda_3 , s,t,r ,l,k\in {\bf N}$,
with $ s \geq 1$,
$\lambda_1 \geq \lambda_2 \geq \lambda_3 $ and maps $SL(V)$-invariant 
and with all the components nonzero:
{\small
$$ 0 \rightarrow S^{\lambda_1, \lambda_2, \lambda_3 } V (t) \rightarrow 
S^{\lambda_1+s, \lambda_2, \lambda_3 } V (t+s) 
\rightarrow E \rightarrow 0 $$
$$ 0 \rightarrow 
S^{\lambda_2 +s-1, \lambda_2, 
\lambda_3 } V (t +\lambda_2 +s-1- \lambda_1)
 \rightarrow 
S^{\lambda_1, \lambda_2, \lambda_3 } V (t) \rightarrow 
S^{\lambda_1, \lambda_2 +s, \lambda_3 } V (t+s)
\rightarrow E \rightarrow 0 $$
$$ 0 
\rightarrow 
S^{\lambda_1-l, \lambda_2-k, \lambda_3 } V (t-k-l) 
\rightarrow 
S^{\lambda_1, \lambda_2-k, \lambda_3 } V (t-k) 
\oplus 
S^{\lambda_1-l, \lambda_2, \lambda_3 } V (t-l) 
\rightarrow S^{\lambda_1, \lambda_2 , \lambda_3 }
 V (t) \rightarrow E \rightarrow 0 $$}
\hspace*{-0.27cm} Then $E$  is stable.
\end{theorem}

\begin{theorem} \label{stairstablebis}
Let $E$ be a homogeneous bundle on ${\bf P}^3={\bf P}(V) $ whose 
minimal free resolution is one of the following
{\small $$ 0 
\rightarrow 
 \oplus_{i=r,...,1+\epsilon}
S^{\lambda_1 +i, \lambda_2-i, \lambda_3 } V  \stackrel{\psi}{\rightarrow} \oplus_{i=r,...,1}
S^{\lambda_1+i+1, \lambda_2 -i, \lambda_3 }
 V (1)
\rightarrow 
E
\rightarrow 
0 $$
 $$ 0 
\rightarrow 
S^{\lambda_1 +1 -k, \lambda_2 -r -k, \lambda_3 } V (1-r-2k) 
\stackrel{\varphi}{\rightarrow} 
 S^{\lambda_1+r -k, \lambda_2 -r-k, \lambda_3 } V (-2k) 
\oplus \oplus_{i=r,...,1+\epsilon}
S^{\lambda_1 +i, \lambda_2-i, \lambda_3 } V  \stackrel{\psi}{\rightarrow}
$$ $$\rightarrow \oplus_{i=r,...,1}
S^{\lambda_1+i+1, \lambda_2 -i, \lambda_3 }
 V (1)
\rightarrow 
E
\rightarrow 
0 $$}
\hspace*{-0.27cm}
for some $\lambda_1, \lambda_2, \lambda_3 , r
 \in {\bf N}$, $\epsilon
\in \{0,1\}$,
$\lambda_1 \geq \lambda_2 \geq \lambda_3 >0$,  $SL(V)$-invariant  maps,
the only nonzero component of $\varphi$ is the first,
  $\psi|_{S^{\lambda_1 +i, \lambda_2-i, \lambda_3 } V}$
 has only the components  into 
$S^{\lambda_1+i+1, \lambda_2 -i, \lambda_3 } V (1)$
 and $S^{\lambda_1 +i, \lambda_2 -i +1, \lambda_3 } V (1)$
 nonzero and  $\psi|_{S^{\lambda_1+r  -k, \lambda_2 -r-k, \lambda_3 } 
V (-2k)}$ has all the components nonzero.

Then $E$  is stable.
\end {theorem}

Finally, by using again the simplicity of classical staircases,  we prove

\begin{theorem} \label{kerscalasempl}
Let $E$ be a homogeneous bundle on ${\bf P}^3$ 
such that  there 
exist $ \lambda_1,\lambda_2,\lambda_3, c,d \in {\bf N} $ 
with $c \neq d$, $\lambda_1 \geq \lambda_2 \geq \lambda_3$
such that   the  minimal free resolution of $E$ is
$$ 0 \rightarrow S^{\lambda_1,\lambda_2,\lambda_3 } V
 \stackrel{\alpha}{\rightarrow} \oplus_{
s \in K } 
S^{\lambda_1+s_1,\lambda_2+s_2,\lambda_3+s_3, s_4 } V (s_1+s_2+s_3 +s_4)
\rightarrow   E \rightarrow 0 $$
where 
$K= \{s \in {\bf N}^4 |
\; 
s_i \leq \lambda_{i-1} -\lambda_i \; for\; i=2,3 ,\;
s_4 \leq \lambda_{3},\;  s_1+s_2+s_3+s_4=c \; and \;
s_1\neq d \; or \; s_4 \neq 0\}$.

Then $E$ is simple.
\end{theorem}

It seems difficult to  generalize these results to ${\bf P}^n$ since it 
is not clear, at least for me,
 how    to generalize the calculations on the slope in $\S4$. 

The outline of the paper is the following: in \S2 we collect some recalls 
on representation theory and  
the quiver associated to homogeneous bundles on ${\bf P}^3$, in \S3  we make 
some calculations about slopes, which are useful to prove, in \S4, 
that homogeneous bundles whose quiver support 
is a parallelepiped or a classical staircase are stable; in \S5 we study 
the resolutions of such bundles
and finally in \S6 we prove Theorem \ref{kerscalasempl}.

\section{Notation and  recalls}

We recall some facts  on representation theory (see for instance
 \cite{F-H}) and on quivers.

Let $d$ be a natural number and let $ \lambda =
(\lambda_1,...,\lambda_k)$ be a partition of $d$ with $\lambda_1 \geq ...
\geq \lambda_k$.

For any $V$ complex vector space of dimension $n$,
 $S^{\lambda} V$ will denote the  {\bf Schur representation} 
($SL(V)$-representation) associated to $\lambda$
(see Lecture 6 in \cite{F-H}).

The $S^{\lambda} V $ are irreducible  $SL(V)$-representations and it is
well-known that all the irreducible $SL(V)$-representations 
are of this form.

\bigskip

We recall that {\bf Pieri's formula} says that, 
 if $\lambda=(\lambda_1, \lambda_2,...) $ is a partition of a natural
number  $d$ with $\lambda_1 \geq \lambda_2 \geq ..$
and $t$ is a natural number, then $$S^{\lambda} V \otimes 
S^t V = \oplus_{\nu } S^{\nu} V$$  as $SL(V)$-representation,
where the sum is performed on all the 
partitions $\nu= (\nu_1,..) $ with $\nu_1 \geq \nu_2\geq ...$ 
 of $d+t$ whose Young diagrams are
obtained from the Young diagram of $\lambda$ adding  $t$ boxes
not two in the same column.

Observe  that, 
if $V$ is a complex vector space of
dimension $n$, then $S^{(\lambda_1,...,\lambda_{n-1})} V$
  is isomorphic, as $SL(V)$-representation,
 to $S^{(\lambda_1+r,...,\lambda_{n-1} +r,r) }V$  for all $r \in 
{\bf N}$. 
Besides
$(S^{(\lambda_1,...,\lambda_{n})}V)^{\vee}$ is isomorphic, as
$SL(V)$-representation, to $S^{ (\lambda_1 -\lambda_n,..., 
\lambda_1-\lambda_{2}) }V$.

\begin{rema} 
a) If $U,W,V$ are three  vector spaces, then 
$Hom (U \otimes {\cal O}_{{\bf P}(V)} (-s)  , W  \otimes {\cal O}) 
= Hom (U \otimes S^s V   ,W) $ (the 
isomorphism is given by $H^0 (\cdot^{\vee})^{\vee} $).

b) Let $V$ be a vector space. For any $\lambda$, $\mu $ partitions, 
$s \in {\bf N}$, up to multiples there is a unique
 $SL(V)$-invariant map  
$$ S^{\lambda} V \otimes {\cal O}(-s) \rightarrow S^{\mu} V \otimes {\cal O}$$
by part a of the remark, Pieri's formula and  Schur's lemma. 
\end{rema}

\begin{lemma} \label{iniet} Let $V$ be a complex vector space of
 dimension 
$n$. Let $\lambda_1,...., \lambda_{n-1}, s \in {\bf N} $ with 
$ \lambda_1 \geq ...  \geq \lambda_{n-1} $. 
 On ${\bf P}^{n-1} ={\bf P}(V) $  any $SL(V)$-invariant
nonzero map $$ S^{\lambda_1,..., \lambda_{n-1}} 
V (-s) \rightarrow S^{\lambda_1 +s , \lambda_2 ,..., \lambda_{n-1}} V$$ 
 is injective.
\end{lemma}

The above lemma is well known; for the proof see for instance \cite{O-R1}.

\begin{notat}
$\bullet $ 
In all the paper  $V$ will be a complex vector space of dimension $4$ if 
not otherwise specified.

$\bullet$ If $E$ is a vector bundle on ${\bf P}(V) $ then $\mu(E)$ will 
denote the slope of $E$, i.e. the first Chern class divided by the rank.
\end{notat}

\begin{defin} (See  \cite{Sim}, \cite{King}, \cite{Hil1}, \cite{G-R}.)
A {\bf quiver} is an oriented  graph ${\cal Q}$ with the set 
${\cal Q}_0$ of 
vertices (or points) and the set  ${\cal Q}_1$ of arrows.

A {\bf path} in ${\cal Q}$ is  a formal composition 
of arrows $ \beta_m ...\beta_1$ where the source of an arrow $\beta_i $
is the sink of the previous arrow $\beta_{i-1}$. A {\bf cycle} is a path 
whose
source of the first arrow $\beta_1$ is the sink of the last arrow 
$ \beta_m$.

A {\bf relation} in ${\cal Q}$ is a linear form $\lambda_1 c_1+...+
\lambda_r c_r$ where $c_i$ are paths in ${\cal Q}$ with a common 
source and
 a common sink and $\lambda_i \in {\bf C}$.

A {\bf representation of a quiver} ${\cal Q} =({\cal Q}_0, 
{\cal Q}_1)$, or {\bf ${\cal Q}$-representation},
  is the couple of a set of vector 
spaces $\{X_i\}_{i \in {\cal Q}_0} $ and of a set of linear maps 
$ \{\varphi_{\beta} \}_{\beta \in {\cal Q}_1}$ where  $\varphi_{\beta} : X_i 
\rightarrow X_j$ if $\beta$ is an arrow from $i$ to $j$.

A {\bf representation of a  quiver ${\cal Q}$ with relations ${\cal R}$} 
is a ${\cal Q}$-representation such that  
$$\sum_j \lambda_j \varphi_{{\beta}^j_{m_j}} ...\varphi_{{\beta}^j_{1}}=0$$ 
for every $ \sum_j \lambda_j \beta^j_{m_j} ...\beta^j_{1} \in {\cal R}$.

Let $(X_i, \varphi_{\beta})_{i \in {\cal Q}_0,\; \beta \in {\cal Q}_1}$ and 
$(Y_i, \psi_{\beta})_{i \in {\cal Q}_0,\; \beta \in {\cal Q}_1}$ be two 
representations of a quiver ${\cal Q}= ({\cal Q}_0 , {\cal Q}_1)$.
A {\bf morphism} $f $ from $(X_i, \varphi_{\beta})_{i \in {\cal Q}_0,\; 
\beta \in {\cal Q}_1}$ to
$(Y_i, \psi_{\beta})_{i \in {\cal Q}_0,\; \beta \in {\cal Q}_1}$ is a
 set of linear maps $f_i : X_i \rightarrow Y_i $, 
$i \in {\cal Q}_0$ such that , for
 every  $\beta \in {\cal Q}_1$,  $\beta$  arrow from $ i$ to 
 $j$, we have $ f_j \circ \varphi_{\beta} = \psi_{\beta} \circ f_i$.

A morphism $f$ is injective if the $f_i$ are injective.
\end{defin}

\begin{notat} \label{multord}
We  say that a representation $(X_i, \varphi_{\beta})_{i \in 
{\cal Q}_0,\; \beta \in {\cal Q}_1}$ 
of a quiver ${\cal Q}= ({\cal Q}_0 , {\cal Q}_1 )$ has 
{\bf multiplicty}  $m$  in a point $i$  of ${\cal Q}$  if $dim \;
X_i =m$.
 
The  {\bf support} (with multiplicities)
 of a representation of a quiver ${\cal Q}$
is the subgraph
of ${\cal Q}$  constituted by the points  of multiplicity  $\geq 1$ and the
 nonzero arrows (with the multiplicities
 associated to every point of the subgraph).

If the quiver has no cycles,
we introduce the following {\bf partial order} on ${\cal Q}_0$:
we say that $A > B$ if there is a path from $B$ to $A$ (that is whose 
source is $B$ and whose sink is $A$).
\end{notat}

\medskip

Observe that  $ {\bf P}^n$ can be seen as $  SL(n+1) /P$,
where 
$$ P= \{ A \in SL(n+1)|  \; A_{i,1}=0 \; for \; i=2,..,n+1
\} $$ Since $P$ is parabolic, we have Levi decomposition $P=RN$, where 
$N$ is unipotent and $R$  is reductive. Here 
$$ R=  \left\{ \left(\begin{array}{cc} a & 0  \\
0    & A   \end{array}\right) \in P, \; a \in {\bf C}, \; A\in M( n \times n,
{\bf C})\right\}$$ 
Let $ {\cal P}$, ${\cal R}$ and ${\cal N}$ be the Lie algebras associated to 
$ P,R,N$. 

We recall that the homogeneous bundles on ${\bf P}^n=SL(n+1)/P$
are given by the representations of $P$. 
A representation of $P$ is completely reducible if and only if it is trivial 
on $N$ (see \cite{Ise}). 
In this case the representations are determined by their restrictions 
on $R$. So  
the irreducible homogeneous bundles on   $ {\bf P}^n$ are given 
by the irreducible $R$-representations. One can easily prove that they are 
the bundles obtained by applying the Schur functors to the quotient bundle
 $Q = T_{{\bf P}^{n} }(-1)$ (where $ T_{{\bf P}^{n}}$
 is the holomorphic tangent bundle) and twisting, that is the bundles
$S^{l_1,...,l_{n-1}} Q (t) $, for 
 $ l_1,...,l_{n-1} \in {\bf N}$, $l_1 \geq ... \geq l_{n-1}$, 
$ t \in {\bf Z}$
(in fact the representation $R \rightarrow gl(1) $ sending 
$ \left(\begin{array}{cc} a & 0  \\
0    & A  \end{array} \right)
$ to $a I$  gives the bundle $ {\cal O} (-1)$ and the representation
$R \rightarrow gl(2) $ sending $\left( \begin{array}{cc} a & 0  \\
0    & A   \end{array} \right)
$ to $A$ gives the bundle $Q$).

Since $Q$ is of rank $n$, the rank of $S^{l_1,...,l_{n-1}} Q(t)$ is the 
dimension  of $S^{l_1,...,l_{n-1}} {\bf C}^n$, therefore
$$rk(S^{l_1,...,l_{n-1}} Q(t))
 =\prod_{1\leq i< j \leq n} \frac{l_i - l_j+j-i}{j-i} \;\;\;\;\;\;\;l_n:=0$$
(see for instance \cite{F-H}, Theorem 6.3). Furthermore we have:
$$\mu(S^{l_1,...,l_{n-1}} Q(t)) =
\frac{l_1+...+l_{n-1}}{n}+t$$ 
in fact: by Euler sequence, $c_1(Q)=1$ and then 
$\mu (Q)= \frac{1}{n}$; besides,
for any vector bundle $E$ of rank $r$, we have  that 
$\mu (S^{\lambda} E) = |\lambda | \mu(E) $ 
for any $ \lambda=(\lambda_1, ...,\lambda_{r-1})$ with  
  $\lambda_1 \geq  ... \geq \lambda_{r-1}$, where $ |\lambda|= 
\lambda_1+ ...+\lambda_{r-1}$
(see  \cite{Ott} Appendix E; one can prove this formula
 for instance by proving  the result first in the case of the 
symmetric powers of $E$, by using the splitting principle,
and then  in general 
by double induction on the number of the rows of the Young diagram and 
the number of the elements in the last row,
by using Pieri's formula  applied to $S^{\lambda} E \otimes S^m E$, for 
 $ \lambda=(\lambda_1, ...,\lambda_{r-1})$, $m \in {\bf N}$).

Finally, from the formulas for the rank and the slope, we get:
 $$c_1(S^{l_1,...,l_{n-1}} Q(t)) =
\left(\prod_{1\leq i< j \leq n} \frac{l_i - l_j+j-i}{j-i} \right)
\left(\frac{l_1+...+l_{n-1}}{n}+t\right)$$

In 1990 Bondal and Kapranov introduced an equivalence between the category of 
homogeneous bundles on certain homogeneous varieties $X$ 
and the category of representations of a certain quiver with relations
associated to the variety $X$, see \cite{B-K}. 
Then, in \cite{Hil1}, \cite{Hil2},  Hille  pointed out that this
 equivalence is not always true with the relations given in \cite{B-K}; 
he introduced a different quiver with relations 
(Hille's quiver has less arrows than Bondal and Kapranov's one, but they 
coincide if $X$ is Hermitian symmetric) 
and proved that the relations are quadratic if the variety is Hermitian 
symmetric; besides he described explicitly the relations for ${\bf P}^2$, 
precisely he showed that, in this case, the relations correspond to the
 commutativity of all square diagrams.
Later, in 2003,  Alvarez-Consul and Garcia-Prada  corrected  the relations on 
Bondal and Kapranov's quiver appropriately and proved the equivalence 
between the category of homogeneous bundles 
and the category of representations of this quiver with relations, 
see \cite{A-G1} Corollary 1.13.
We mention also the later paper   \cite{O-R2} (and in particular Theorem 5.9 
and Corollary 8.5) to which we mainly  refer. In it there is an explict 
description of the relations in the case of Grassmannians 
(see Proposition 8.4); in particular Hille's result on the relations 
in the case of ${\bf P}^2$ is extended
to $ {\bf P}^n$.

Here we state the result  (Theorem \ref{BKH}) only for $X={\bf P}^n$, since 
we need only this case.

\begin{defin} \label{defQ} Let
 ${\cal Q}= {\cal Q}_{{\bf P}^n}$ be defined in the following way:

$\bullet $ let $${\cal Q}_0 = \{irreducible \;R-representations\}=$$
$$\{irreducible \;hom. \; bundles\; on \;
{\bf P}^n\}=  \{S^{l_1,...,l_{n-1}} Q (t) 
| \; l_1,...,l_{n-1} \in {\bf N},\; l_1 \geq ... \geq l_{n-1}, 
\; t \in {\bf Z}\} $$

$\bullet $
let ${\cal Q}_1$ be defined in the following way: there is an arrow from 
$ S^{l_1,...,l_{n-1}} Q (t) $  to $  S^{l'_1,...,l'_{n-1}} Q (t')  $
 iff 
$  S^{l'_1,...,l'_{n-1}} Q (t')  $,
is a direct summand of  
\begin{center}
$ \Omega^1 \otimes S^{l_1,...,l_{n-1}} Q (t)=
 \wedge^{n-1} Q (-2) \otimes S^{l_1,...,l_{n-1}} Q (t)= $
$=S^{l_1+1,...,l_{n-1}+1} Q(t-2) \oplus \oplus_{i=1,...,n-1}
S^{l_1+1,...,l_{i-1}+1,l_i,l_{i+1}+1,
...,l_{n-1}+1,1} Q( t-2) = $
$=S^{l_1+1,...,l_{n-1}+1} Q(t-2) \oplus \oplus_{i=1,...,n-1}
S^{l_1,...,l_{i-1},l_i -1,l_{i+1},
...,l_{n-1}} Q( t-1)  $
\end{center}

\end{defin}

(The last equality holds because, 
by the Euler sequence, $\wedge^{n} Q = {\cal O}(1)$).

Our quiver has $n+1$ connected components ${\cal Q}^{(1)}$,  
....,  ${\cal Q}^{(n+1)}$  (see Theorem 5.12 in \cite{O-R2}); they are
given by the congruence class modulo 
$(n+1)/n (= -\mu(\Omega^1)) $ of the slope of the homogeneous
 vector bundles corresponding to the points of the connected component.

According to the definition of our quiver, from 
the point of the quiver corresponding to 
$S^{l_1,...,l_{n-1}} Q (t)$ we can have $n$ arrows, one toward 
 $S^{l_1+1,...,l_{n-1}+1} Q(t-2)$ and, for any 
$i =1,..., n-1$, one toward 
$S^{l_1,...,l_{i-1},l_i -1,l_{i+1},...,l_{n-1}} Q( t-1)$. 
So we can identify the points of every connected component ${\cal Q}^{(r)}$  
of ${\cal Q}$    with a subset of ${\bf Z}^{n}$. We call
$\{V_i \}_{i=0,..., n}$ its canonical basis, precisely, for 
any $i=1,..., n-1$, we call $V_i$  
 the ``not applied'' vector corresponding to the arrows of the quiver
 from the points  
 $ S^{l_1,...,l_{n-1}} Q (t) $
 toward the points  $S^{l_1,...,l_{i-1},l_i -1,l_{i+1},...,l_{n-1}} Q( t-1) $
and we call $V_0 $  the vector corresponding to the arrows 
 from the points $ S^{l_1,...,l_{n-1}} Q (t) $ toward the points 
 $S^{l_1+1,...,l_{n-1}+1} Q(t-2)$.

The figure on the left  shows one of the three connected
 components in the case $n=2$. 
The figure on the right shows the border planes of 
one of the four connected components in the case $n=3$ ($t=0,1,2,3$); they are the
planes $ \sigma= {\cal O}(t)+\langle V_0 , V_1+V_2\rangle=\{S^{k,k}Q(s) | \;
\mu (S^{k,k} Q(s)) \equiv \mu ({\cal O}(t)) \;mod\; (n+1)/n  \}$ 
and
 $ \pi= {\cal O}(t)+\langle V_1 , V_0+V_2\rangle=
\{S^{k}Q(s) | \;
\mu (S^{k} Q(s)) \equiv \mu ({\cal O}(t))\; mod \; (n+1)/n \}$.

\begin{center}
\includegraphics[scale=0.3]{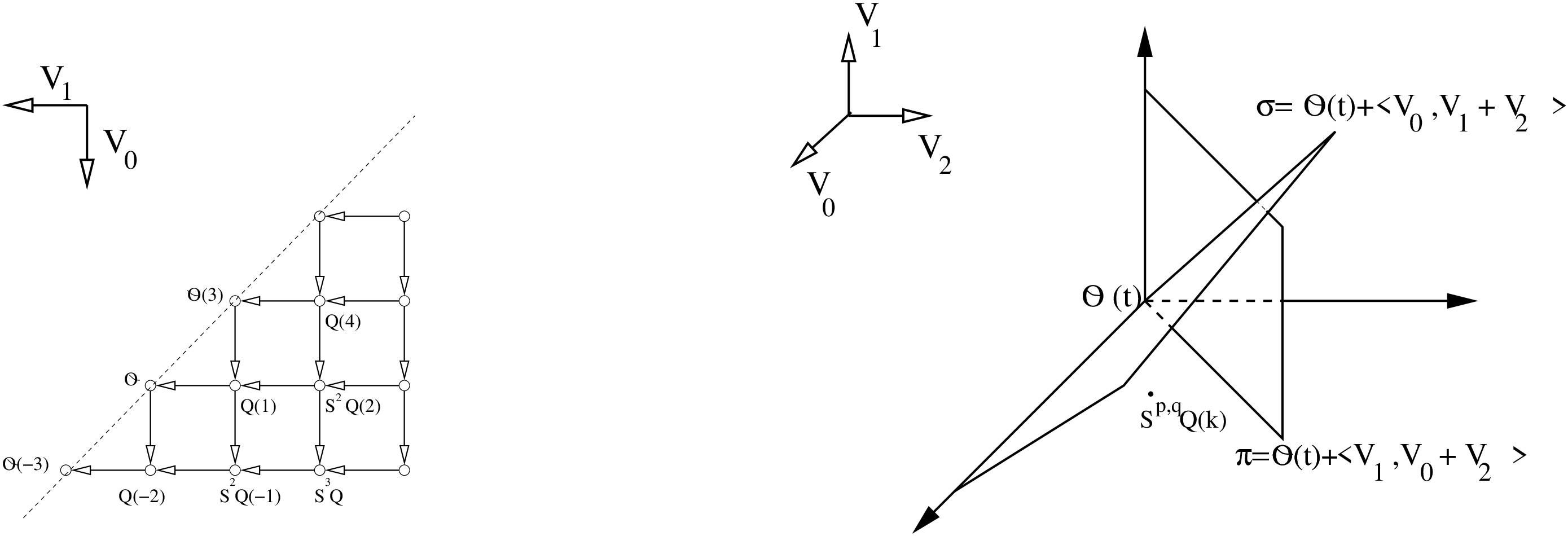}
\end{center}

\begin{defin} \label{defR} 
Let ${\cal R}$ be the relations on ${\cal Q}={\cal Q}_{{\bf P}^n}$ 
(defined in Definition \ref{defQ})
 given by the commutativity of  the squares, that is, for any $ P   \in {\cal Q}_0^j  $ for some $j$,  
$$ \beta_{ P+V_i+V_k , P +V_k }
\beta_{  P+V_k , P } - \beta_{ P+V_i+V_k , P +V_i }
\beta_{  P+V_i , P }$$ 
 and 
 $$\beta_{ P+V_i+V_k , P +V_k }
\beta_{  P+V_k , P }$$ if
  $ P + V_i \not \in {\cal Q}_0$ (where $\beta_{T,P}$ denotes
 the arrow from $P$ to $T$). 
\end{defin}

For every $E$ homogeneous bundle on  ${\bf P}^n$, 
let  $gr E$ be the bundle given by the  restriction to $R$
of the representation of   $P$ which gives  $E$.

A homogeneous bundle $E$ on ${\bf P}^n$ gives a ${\cal P}$-representation,
 which we call again  $E$, and
given a ${\cal P}$-representation $E$, the action of  ${\cal N}$
on $E$ induces a morphism of ${\cal R}$-representations
\begin{center}$\theta_E :{\cal N}\otimes gr E \rightarrow gr E$
\end{center}
\begin{defin}  (See Definition 5.5 in \cite{O-R2}.)
Let ${\cal Q}$ be the quiver defined in Definition \ref{defQ}.
 We define the ${\cal Q}$-representation 
 associated to a homogeneous bundle $E$
 on $ {\bf P}^n$ in the following way.
Let $grE=\oplus_{F \in {\cal Q}_0} F \otimes V_{F}$.

$\bullet $ 
To $F \in{\cal Q}_0 $ we associate the vector space   $V_{F}$.

$\bullet $ 
To fix the linear maps associated to the arrows: consider the 
decomposition
$$Hom (gr E \otimes\Omega^1, gr E)=\oplus_{F, F'}
 Hom(V_{F},V_{F'})
\otimes Hom (F \otimes \Omega^1,F')$$
We fix  generators $m_{F'F}$ of 
$Hom(F\otimes\Omega^1, F')^G$:
fix  a maximal vector 
$v_{F}\in F$ $\forall F\in{\cal Q}_0$;
let $ \xi_i$ $i=1,...,n$  be  the weights of
 the ${\cal R}$-representation
 ${\cal N}$ and $ n_i\in{\cal N}$ corresponding maximal 
vectors;
we can prove that $Hom(F\otimes\Omega^1, F')^G
$ is equal to $0$ or ${\bf C}$ and it is equal to  ${\bf C}$
iff the difference of the maximal weights of $ F$ and  $ F'$ is $ 
  \xi_i$  for some $i$; 
fix  a generator $m_{F'F}$ of
$Hom(F\otimes\Omega^1, F')^G$
taking  $v_{F}\otimes n_j$ to $v_{F'}$.

Then
$\theta_E \in Hom (gr E \otimes\Omega^1, gr E)^G$ can be written
 as $$ \sum_{F, F'} g_{ F', F}
\otimes m_{F', F}$$
$ g_{F', F } \in Hom(V_{F}, V_{F'})$ 
is the linear map associated to the arrow from $ F $ to $ F'$ 

\end{defin}

As we have already said,  many people contributed to the following theorem. 
The explicit description of the relations for $ {\bf P}^n$  
can be found in  \cite{O-R2}, see Corollary 8.5.

\begin{theorem} \label{BKH}
The category of the homogeneous bundles on $  {\bf P}^n$ is equivalent
to the category of finite dimensional representations of the 
quiver ${\cal Q}$ defined in Definition \ref{defQ} with the relations 
${\cal R}$ defined in Definition \ref{defR}.
\end{theorem}

\begin{notat} $\bullet $ 
We will often speak of the {\bf ${\cal Q}$-support} 
of a homogeneous bundle $E$ instead of the  support with multiplicities 
of the ${\cal Q}$-representation of $E$.

 $\bullet $
The word ``{\bf parallelepiped}'' will denote the subgraph with multiplicities 
of ${\cal Q}$ given by the subgraph of ${\cal Q}$  included
 in a parallelepiped whose sides are parallel to $\langle
V_{i_1},..., V_{i_{n-1}} \rangle$ for some distinct $i_1,...,i_{n-1}$, with the
multiplicities  of all its points equal to $1$.

 $\bullet $ If $A$ and $B$ are two subgraphs of ${\cal Q} $, $A \cap B$ is the
 subgraph of ${\cal Q}$ whose vertices and arrows are the vertices and arrows
both of $A$ and of $B$; $A -B$ is the subgraph of ${\cal Q} $ whose vertices 
are the vertices of $A$ not in $B$ 
   and the arrows are the arrows of $A$ joining two vertices of $A -B$.
\end{notat}

\begin{rema} \label{Spq}  (\cite{B-K}) 
The ${\cal Q}$-support of $S^{p_1,...,p_n}V(t)$ is  
a parallelepiped with vertex with maximum slope
$S^{p_1 -p_n  ,...., p_{n-1} -p_n} Q (t +p_n )$
and the side of direction $V_1$ of length $p_1 -p_2 $,.....,
the side of direction $V_{n-1}$ of length  $p_{n-1} -p_n$,
the side of direction $V_0$ of length  $ p_n$.  

If $ n=3$,  for the parallelepiped support of $ S^{l_1,l_2,l_3 } V(t)$, 
 the vertex opposite to the vertex  with maximum slope
 in the side in the direction 
$ \langle V_1, V_2 \rangle $ corresponds to 
$ S^{l_2 -l_3 } Q (-l_1 +2l_3 +t)$ (see the figure below).
\end{rema}

In fact:  by  the Euler sequence $ S^{p_1,...,p_n} V (t)
= S^{p_1,...,p_n} ({\cal O} (-1) \oplus Q )(t) $ as $R$-representation; 
by  the formula of a Schur functor applied to a direct sum (see
\cite{F-H}, Exercise 6.11) we get 
$$S^{p_1,...,p_n} V = \oplus S^{\lambda} Q \otimes S^m {\cal O}
(-1)   $$ as $R$-representations, 
 where the sum is performed on $ m \in {\bf N}$ and on $ \lambda $ Young 
diagram obtained  from the Young diagram of $(p_1,...,p_n)$ by
 taking off $m$ boxes not
 two in the same column; thus 
 $$
S^{p_1,...,p_n} V  =\oplus_{ 
{\tiny \begin{array}{c}
0 \leq m_1 \leq p_1 -p_2  ,\\
.....,\\
0 \leq m_{n-1} \leq p_{n-1} -p_n, \\
 0 \leq m_n \leq p_n  
\end{array}}} 
S^{p_1 -m_1,...., p_n-m_n} Q (-m_1 -....-m_n +t)= 
$$ 
$$=\oplus_{{\tiny \begin{array}{c}
0 \leq m_1 \leq p_1 -p_2  ,\\
.....,\\
0 \leq m_{n-1} \leq p_{n-1} -p_n, \\
 0 \leq m_n \leq p_n  
\end{array}
}}
S^{p_1 -p_n -(m_1-m_n) ,...., p_{n-1} -p_n-(m_{n-1}-m_n)} 
Q (-m_1 -....-m_n +t +p_n -m_n)
$$ 
Finally to show the maps associated to the arrows in the 
parallelepiped are nonzero 
we can argue as in the proof of  Remark 23 in 
\cite{O-R1}

\begin{center}
\includegraphics[scale=0.36]{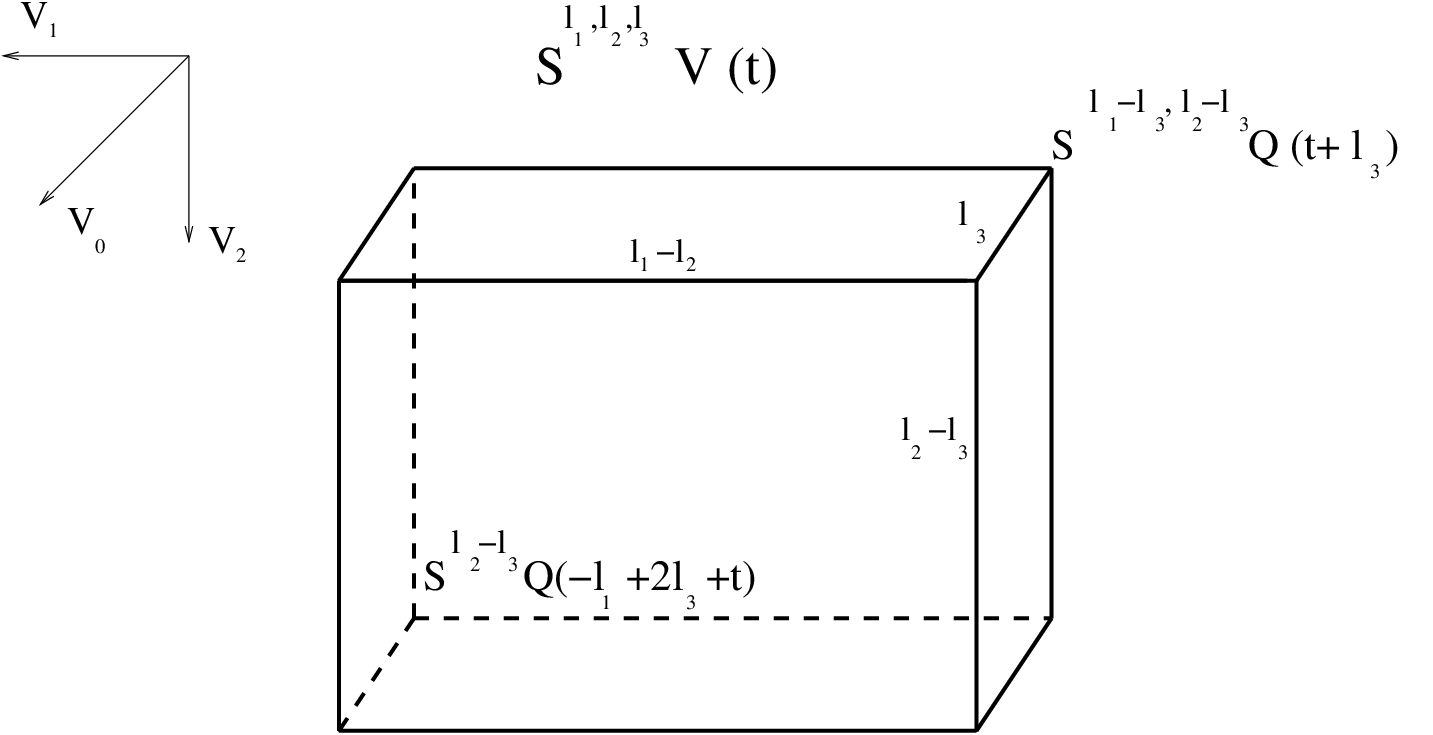}
\end{center}

\section{Some technical lemmas to study stability}

\begin{rema}  \label{gio}  i) The first Chern class
  of a homogeneous bundle $E$ can be calculated as the 
sum of the first Chern classes of the irreducible bundles corresponding to 
the vertices of the ${\cal Q}$-support of $E$ multiplied 
by the multiplicities. The rank of $E$ is  the 
sum of the ranks of the irreducible bundles corresponding to 
such vertices  multiplied by the multiplicities.

We will often speak of the slope (resp. $c_1$, rank) of  a 
graph with multiplicities  instead of
 the slope (resp. $c_1$, rank) of the  vector 
bundle whose ${\cal Q}$-support is  that graph with multiplicities.

ii) Suppose  the set of the vertices of the 
${\cal Q}$-support of $E$ is the disjoint union of the 
vertices of the supports of two ${\cal Q}$-representations $A$ and $B$;
if $\mu(A)= \mu(B) $ then $\mu(E)= \mu(A) = \mu(B)$, if  
 $\mu(A) < \mu(B) $ then $\mu(A) < \mu(E) < \mu(B)$.
In fact, by i, $$ \mu (E)=\frac{c_1(E)}{rk(E)}= 
\frac{c_1(A)+ c_1(B)}{rk(A)+rk(B)}$$
so, if  $\frac{c_1(A)}{rk(A)}<\frac{c_1(B)}{rk(B)}$, we get immediately
$\frac{c_1(A)}{rk(A)} < \mu(E) < \frac{c_1(B)}{rk(B)}$; 
if $\frac{c_1(A)}{rk(A)}=\frac{c_1(B)}{rk(B)}$, there exists 
$t$ such that $c_1(B)= t \; c_1 (A)$ and  $rk(B)= t\; rk(A)$ and 
we conclude easily by substituting.

iii)
On ${\bf P}^n$ we have 
$S^{l_1,...,l_{n-1}} Q(t)^{\vee} = 
S^{l_1,l_1 -l_{n-1},...,l_1-l_2} Q(-t-l_1)$. 

In particular on ${\bf P}^3$ $S^{l_1,l_2} Q(t)^{\vee} = 
S^{l_1,l_1-l_2} Q(-t-l_1)$. Thus, by dualizing, the 
vector $ V_0 $ goes into  $-V_1$ and $V_2 $ into $ -V_2$. 
\end{rema}

\begin{lemma} \label{caldo}
Let $c \in {\bf N}$. Let $E_{S^{l_1,l_2}Q(t),c , V_i , V_k}$ 
be the hypotenuse of the isosceles triangle
 whose vertices are: 
$$A := S^{l_1,l_2} Q (t) \;\;\;\;\;
 A + c V_i \;\;\;\;\;   A+ c V_k$$ 
for some $i,k \in \{0,1,2\}$. 
Let  $R_{S^{l_1, l_2 }Q (t) , c , V_i, V_k}$
 be the hypotenuse of the isosceles
triangle whose vertices are: 
$$A := S^{l_1,l_2} Q (t) \;\;\;\;\; A - c V_i \;\;\;\;\; A- c V_k$$
Let $x = l_1 - l_2 +1$ and $z = l_2 +1$. 

Let $e_{S^{l_1,l_2}Q(t),c , V_i , V_k}$ 
and  $r_{S^{l_1,l_2}Q(t),c , V_i , V_k}$ be the sum 
of the  ranks of the vector bundles 
corresponding to the points of  $E_{S^{l_1,l_2}Q(t),c , V_i , V_k}$ 
and of  $R_{S^{l_1,l_2}Q(t),c , V_i , V_k}$ respectively. 
We have:

a) $ e_{S^{l_1,l_2}Q(t),c , V_1 , V_2} = 
(c+1) [-c x (x+2z +1) +2xz (x+z) ]$

b) $
r_{S^{l_1,l_2}Q(t),c , V_1 , V_2}  = (c+1) [c x(x+2z -1 ) +2zx (z+x)]$

c) $
 e_{S^{l_1,l_2}Q(t),c , V_0 , V_1} = -c^2 + c(x+z)(x-z+1) +2xz (x+z)$
\end{lemma}

{\it Proof.}  It is a freshman's calculation.
\hfill \framebox(7,7)

\begin{rema} \label{semplice}
Let $a,b,c,d , r,s \in {\bf R}$ with $r$ and $s$ positive. Suppose
$\frac{a}{b} > \frac{c}{d}$. Then
$\frac{ra+sc}{rb+sd} > \frac{a+c}{b+d}$ if $s < r$   and 
$\frac{ra+sc}{rb+sd} < \frac{a+c}{b+d}$ if $s > r$.   
\end{rema}

\begin{prop} \label{main} Let $S$ be a  rectangle in  ${\cal Q}^{(r)}$
for some $r$ (one of the  components of the quiver). Suppose the rectangle
is parallel to $\langle V_i , V_k \rangle$ for some $i,k \in 
\{0,1,2\}$, $i \neq k$. 
Let $S'$ be   obtained by translating 
$S$ in  ${\cal Q}^{(r)} \subset {\bf Z}^3$ 
by $V_j - V_i  $ or by $V_j - V_k  $ where $j \in 
\{0,1,2\} $, $j \neq i,k$.   
Then $\mu(S) > \mu(S')$ with only two exceptions:

a) $\{i,k\} = \{1,2\}$, $ S' $   is  obtained by translating 
$S$ by $V_0 - V_1  $ and the length of the side of $S$ with direction 
$ \langle V_2 \rangle$ is greater than
 the length of the side of $S$ with direction 
$ \langle V_1 \rangle$

b) $\{i,k\} = \{1,2\}$, $ S' $   is  obtained by translating 
$S$ by $V_0 - V_2  $ and the length of the side of $S$ with direction 
$ \langle V_1 \rangle$ is greater than
 the length of the side of $S$ with direction 
$ \langle V_2 \rangle$.
\end{prop}

{\it Proof.}
To prove the statement, we will consider the quotient $\nu$ of
 the sum of the ranks of the bundles 
 corresponding the points   of $S'$ with a fixed $\mu$  and 
 the sum of the ranks of the bundles 
 corresponding the points   of $S$ with the same  $\mu$.
We will show that this quotient is a decreasing function of $\mu$. 
By Remark \ref{semplice} this is sufficient to prove our statement. 

Obviously the points of  $S$ corresponding to  bundles 
  with a fixed $\mu$  form a segment forming an angle of $45$ degrees 
with the sides of $S$.

We have to consider three cases: 
1) $S$ is parallel to  $\langle V_1 , V_2 \rangle$, 
2) $S$ is parallel to  $\langle V_0 , V_1 \rangle$, 
3) $S$ is parallel to  $\langle V_0 , V_2 \rangle$. 

\medskip

Observe that if we consider the segment given by 
the  points of $S$  with a certain $\mu$ and the segment given by 
the  points of $S$  with the slope equal to $\mu -\mu ( \Omega^1)$, 
we have the four subsubcases  A,B,C,D shown in the picture:

\begin{center}
\includegraphics[scale=0.39]{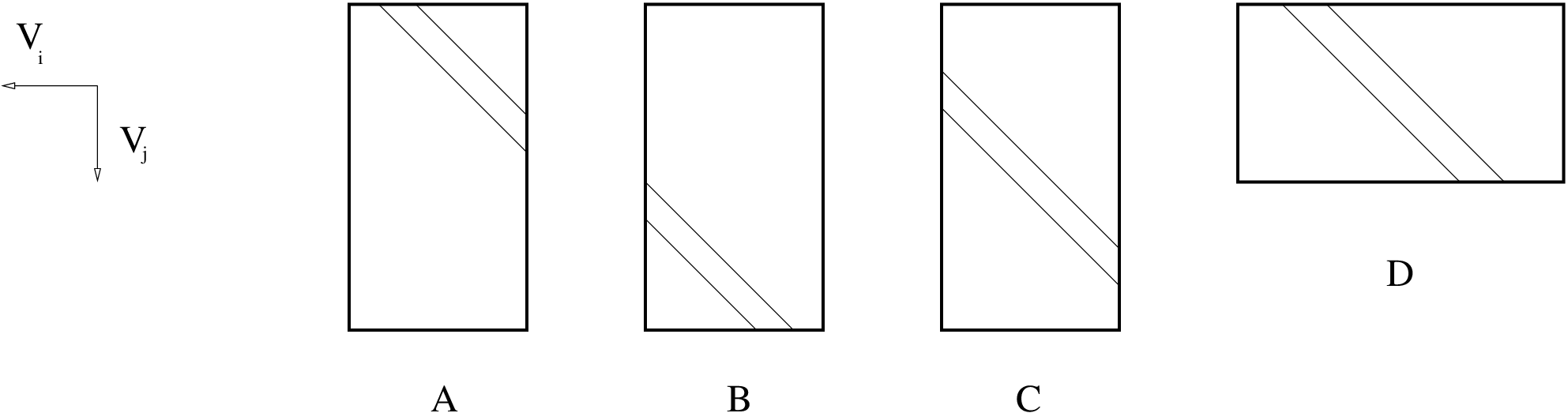}
\end{center}

In cases A,C,D 
let $ S^{l_1,l_2} Q (t)$ be the vector bundle corresponding to the vertex 
of $S$ with greatest slope. In case B let $ S^{l_1,l_2} Q (t)$ be
 the vector bundle corresponding to the vertex 
of $S$ with least slope.
Let $x = l_1 - l_2 +1$, $z = l_2 +1$ and $y=x+z$.

{\bf 1) ${\bf S}$ IS PARALLEL TO  ${\bf
\langle V_1 , V_2 \rangle}$. }

We  divide this case into two main subcases:
the case where $S'$ is   obtained by translating $S$ 
by $V_0 - V_1  $ and the case where   $S'$ is   obtained by translating 
$S$ by $V_0 - V_2  $.

{\bf $\bullet$ Subcase: ${\bf S'}$ is   obtained by translating ${\bf S}$ 
by ${\bf V_0 - V_1 } $.}

\smallskip

{\bf  A)}
Let  $ \nu (c) = 
\frac{
e_{S^{l_1,l_2}Q(t),c , V_1 , V_2}
}{ 
e_{S^{l_1+2,l_2+1}Q(t-1),c , V_1 , V_2} 
}=
\frac{-c (x+1) (x+2z +4) +2(x+1)(z+1) (x+z+2) }{-c x (x+2z +1) +2xz (x+z)}$
(see Lemma \ref{caldo} for notation).

We have to prove that $ \nu(c)$ is an increasing function of $c$.
The numerator  of the derivative $ \nu' (c) $ is
$$ 2 x (x+1)[-3 z(x+z)   + (x+2z+1)  (x+3z+2) ]$$
which is obviously positive.

\smallskip

{\bf  B)}
Let $ \nu (c) = 
\frac{
r_{S^{l_1,l_2}Q(t),c , V_1 , V_2}
}{ 
r_{S^{l_1+2,l_2+1}Q(t-1),c , V_1 , V_2} 
}=
\frac{c  (x+1) [(x+1) +2(z+1) -1 ] +2 (z+1) (x+z+2)}{c
 x (x+2z -1) +2xz (x+z)}$.

We have to prove that $ \nu(c)$ is a decreasing function of $c$.
The numerator of  $ \nu' (c) $ is
$$ 2(x+1)x[ 3 z (x+z) -  (x+3z +2) ( x+2z -1)]$$ 
which is obviosly negative.

\smallskip
(We don't consider subsubcase C because it's the exception a of the
 statement of the theorem.)

\smallskip

{\bf  D)}
Let
$\alpha (c,z,y) = \frac{ 
e_{S^{l_1+2,l_2+1}Q(t-1),c , V_1 , V_2}} 
{e_{S^{l_1,l_2}Q(t),c , V_1 , V_2}}= \frac{c(y-z+1) (y+z+4)
-2 (y-z+1)(z+1)(y+2)}{c(y-z) (y+z+1)  -2 (y-z)zy}$.

We have to show that $\frac{\partial \alpha}{
\partial y}$ is negative.
Let $l$ be the number such that $x=lz$ (thus $y = lz+ z$).
The numerator of $\frac{\partial \alpha}{
\partial y}$ expressed in function of $z,c,l$ is
$$ (-12 z^3 +12 z^3 c -4 c^2 z^2 -12 z^4 +4c z^2 )l^2 +$$
$$-4z[(2z^3 +c^2 z   - c z  -3 z^2 c)+(14 cz -8z^2 -6c^2)+
(8 z^3 c - 4c^2 z^2 -4 z^4 +4c +2c^2)]l +$$
$$+4c -4c^2 z^2 +
22 cz^2 - 8z^2 -4z^4 +8 z^3 c -10 c^2 z +14cz -12z^3 -4c^2$$
The coeffficient of $l^2$ is obviously negative.
The coefficient of $l$ 
is negative since
$2z^3 +c^2 z   - c z  -3 z^2 c$, 
$4z -2c$ and
$6 z^2 +2 c^2 -6cz $
 are obviously positive.
The known term can be seen as the sum of 
$22 c z^2 -10 c^2 z -12 z^3$,
$14 cz -8z^2 -6c^2$ and
$8 z^3 c - 4c^2 z^2 -4 z^4 +4c +2c^2$,
which are negative. 

{\bf $\bullet $ Subcase: ${\bf S'}$ is obtained by translating ${\bf S}$ 
by ${\bf V_0 - V_2 } $}

\smallskip

{\bf  A)} 
Let $ \nu (c) = 
\frac{
e_{S^{l_1,l_2}Q(t),c , V_1 , V_2}}{ 
e_{S^{l_1+1,l_2+2}Q(t-1),c , V_1 , V_2} }=
\frac{-c (x-1) (x+2z +4) +2(x-1)(z+2) 
(x+z+1) }{-c x (x+2z +1) +2xz (x+z)}$.

We have to prove that $ \nu(c)$ is an increasing function of $c$.
The numerator of  $ \nu' (c) $ is
$$ 2 x (x-1)[  -3 z (x+z)+(x+2z+1) (z +2x+2z+2) ]$$
which is obviously positive.

\smallskip

{\bf  B)}
Let $ \nu (c) = 
\frac{
r_{S^{l_1,l_2}Q(t),c , V_1 , V_2}
}{ 
r_{S^{l_1+1,l_2+2}Q(t-1),c , V_1 , V_2} 
}=
\frac{c  (x-1) (x+2z +2)  +2 (x-1)(z+2) (x+z+1)}{c
 x(x +2z -1) +2xz (x+z)}$.

We have to prove that $ \nu(c)$ is a decreasing function of $c$.
The numerator of   $ \nu' (c) $ is
$$ 2(x-1)x[( 3 z (x+z)  -(2x+3z +2) ( x+2z -1)]$$  
which is obviosly negative.

\smallskip

{\bf  C)}
Let $\alpha (c,z,y) = \frac{ e_{S^{l_1+1,l_2+2}Q(t-1),c , V_1 , V_2}} 
{e_{S^{l_1,l_2}Q(t),c , V_1 , V_2}}= \frac{c(y-z-1) (y+z+4)
-2 (y-z-1)(z+2)(y+1)}{c(y-z) (y+z+1)  -2 (y-z)zy}$.

We have to show that $\frac{\partial \alpha}{
\partial z}$ is negative.
Let $l$ be the number such that $x=lz$ (thus $y = lz+ z$).
The numerator of $\frac{\partial \alpha}{
\partial y}$ expressed in function of $z,c,l$ is
$$ -8 z^4 l^4 + (-16 z^4 +8z^3 c )l^3
+(-12 z^4 + 8 z^2 +12 z^3 c -4c^2 z^2  -2 cz^2 ) l^2$$
$$+(-8 z^4 -12 z^3 + 12 c z^3   -4 c^2 z^2 +16 cz^2  -2c^2 z  + 6cz)l$$
$$+4c -4c^2 z^2 +22 cz^2 -8 z^2 -4z^4 +8 z^3 c-10 c^2 z+14 cz -12 z^3 -4 c^2
 $$
The coefficients of $l^4 $, $l^3$, $l$ and the known term
are obviously negative (observe that  for $z=1$, 
$c$ must be $0$).
The coefficient of $l^2 $ is negative too, in fact it can be seen as 
the sum of
$-12 z^4 +12 z^3 c$ 
and
$ (-4c^2 -2c +8) z^2$;
the first term is obviously negative, the second is negative for $c \geq 2$.
In the case $c=1$, observe that the coefficient of $l^2 $ becomes
 $-12 z^4 +12 z^3 +2 z^2$ which is negative for $z \geq 2$ (which 
is our case  since $ c=1$ and $ c \leq z-1$).



\medskip

{\bf 2) ${\bf S}$ IS PARALLEL TO ${\bf \langle V_0 , V_1 \rangle}$. }

\smallskip

{\bf $\bullet $ Subcase: ${\bf S'}$ is   obtained by translating ${\bf S}$ 
by ${\bf V_2 - V_0 } $}

\smallskip

{\bf  A)} 
Let $ \nu (c) = 
\frac{
e_{S^{l_1,l_2}Q(t),c , V_0 , V_1}
}{ 
e_{S^{l_1-1,l_2-2}Q(t+1),c , V_0 , V_1} 
}=
\frac{-c^2 +c
 (x+z-1) (x-z +4) +2(x+1)(z-2) (x+z-1) }{-c^2 +c (x+z) (x-z +1) +2xz (x+z)}$.

We have to prove that $ \nu(c)$ is an increasing function of $c$.
The numerator of $\nu'(c)$ is
$$ c^2 (2x +4z -4)+c
 (-4xz -4x^2 +2 z^2 -6z +4)+
 +2 (x+z) (x+z -1) [2x^2 +4x +z^2 -3z +2]$$
which
is obviously positive (observe  that $c$ must be less or equal than $z$).

\smallskip

{\bf  C)} 
Let $\alpha (c,z,y) = \frac{ 
e_{S^{l_1-1,l_2-2}Q(t-1),c , V_1 , V_2}} 
{e_{S^{l_1,l_2}Q(t),c , V_1 , V_2}}= \frac{
-c^2 + c (y-1) (y-z+1) (y-2z+4)  +2 (y-z+1)(z-2)(y-1)}{
-c^2 + c y (y-z) (y-2z+1)  +2 (y-z)zy}$.

We have to show that $\frac{\partial \alpha}{
\partial y}$ is negative.
Let $l$ be the number such that $x=lz$ (thus $y = lz+ z$).
The numerator of $\frac{\partial \alpha}{
\partial y}$ expressed in function of $z,c,l$ is
$$ (-4 z^2 c - 2 c^2 z^2 ) l^2 
+ [ - 8 z^4 + (24 -12c)  z^3  + (- 8c^2   +20c   -16) z^2+( - 8 c
+ 16 c^2) z
] l 
$$ $$
- 4 z^4 + (12 - 4  c)z^3
 +( -22 c^2 + 2c - 8) z^2   + (6c  
 + 2 c^2) z -4 c  -2 c^3  + 4 c^2 
$$ 
The known term is obviously negative if $z \geq 2$ and $c \geq 1$. 
One can easily check that also in the other cases, i.e. $z = 2$, $c=0  $
and $ z = 1$, $c=0 $ (we recall that $ c \leq z-1$), it is negative.
The coefficient of $l^2$ is always less or equal than $0$.
 The coefficient of $l$  is obviously negative if $ c \geq 2$. If $ c=1$ 
then $z $ must be greater or equal than $2$ and also in this case we 
can easily check  that the coefficient of $l$ is negative.

\smallskip

{\bf B)},
{\bf  D)} 
They follow from subcase A and 
C by duality (see Remark \ref{gio} iv).

\smallskip

{\bf $\bullet $ Subcase: ${\bf S'}$ is   obtained by translating 
${\bf S}$ 
by ${\bf V_2 - V_1 } $}

\smallskip

{\bf  A)} 
Let $ \nu (c) = \frac{
e_{S^{l_1,l_2}Q(t),c , V_0 , V_1}
}{ 
e_{S^{l_1+1,l_2-1}Q(t+1),c , V_0 , V_1} 
}=
\frac{-c^2 +c
 (x+z+1) (x-z +4) +2(x+2)(z-1) (x+z+1) }{-c^2 +c (x+z) (x-z +1) +2xz (x+z)}$.

We have to prove that $ \nu(c)$ is an increasing function of $c$.
The numerator of  $ \nu' (c) $
 is
$$c^2
(4x +2z +4)+ c
( -2x^2 - 6x +4xz +4 z^2 -4)+
 (x+z) (x+z +1) [x^2 +3x +2z^2 -4z +2]$$
which
is obviously positive (observe  that $c$ must be less or equal than $z$).

\smallskip

{\bf  C)} 
Let
$\alpha (c,z,y) = \frac{ 
e_{S^{l_1+1,l_2-1}Q(t-1),c , V_1 , V_2}} 
{e_{S^{l_1,l_2}Q(t),c , V_1 , V_2}}= \frac{
-c^2 + c (y+1) (y-z+2) (y-2z+4)  +2 (y-z+2)(z-1)(y+1)}{
-c^2 + c y (y-z) (y-2z+1)  +2 (y-z)zy}$.

We have to show that $\frac{\partial \alpha}{
\partial y}$ is negative.
Let $l$ be the number such that $x=lz$ (thus $y = lz+ z$).
The numerator of $\frac{\partial \alpha}{
\partial y}$ expressed in function of $z,c,l$ is
$$ [-12 z^4 +(-12  c + 12) z^3  + (4c  - 4c^2) z^2 ] l^2+
[ - 16 z^4 -12 z^3 c + (16  - 4 c^2  - 20 c)z^2  +(8c -4 c^2) z ] l$$
$$
- 8 z^4   - 8 z^3 c + (8 - 2c -4 c^2) z^2  +(-2c^2 - 6c)z  - 4 c^3+ 2c^2
+4c 
$$ 
One can easily see that 
every coefficient of the above polynomial in $l$ is negative.

\smallskip

{\bf B}, {\bf  D)} 
They  follow respectively from subcase A and C by duality 
(see Remark \ref{gio} iv).

\medskip

{\bf 3) ${\bf S}$ IS PARALLEL TO  ${\bf \langle V_0 , V_2 \rangle}$. }

It follows from case 1 by duality (see Remark \ref{gio} iv).
\hfill  \framebox(7,7)

\begin{cor} \label{segm}
Let $S$ be a segment and let $S'$ be obtained by translating $S$ by 
$V_i $ with $ i  \in \{1,2,3\}$: $ S' = S + V_i$.  Then 
$\mu (S') < \mu (S)$.  
\end{cor}

{\it Proof.} {\it i)}
If the direction of $S$ is $\langle V_i \rangle$, it is obvious. 

ii) Otherwise  let the direction of  $S$ be  $\langle V_j \rangle$. From
Proposition \ref{main} and part i) respectively we have 
$ \mu (S) > \mu (S +V_i - V_j ) > \mu (S + V_i ) $. 
\hfill  \framebox(7,7)

\begin{cor} \label{rettangoli}
Let $U$ be a rectangle and let $U'$ be obtained by translating $U$ by 
$V_i $ with $ i  \in \{1,2,3\}$.  Then 
$\mu (U') < \mu (U)$.  
\end{cor}

{\it Proof.} {\it i)} When $V_i $ is  contained in the direction of  $ U$ 
we get the statement from Corollary \ref{segm}. 

ii) Suppose  $V_i $ is not  contained in the direction of  $ U$.
Let the direction of $ U$ be $ \langle V_j ,V_k \rangle$ with $i \neq j,k$.
and let's suppose  the length of the side of $U$ with direction 
 $\langle V_j \rangle$ greater than  the length of 
the side of $U$ with direction 
 $\langle V_k \rangle$.
By Proposition \ref{main} 
and part i) respectively we have 
$ \mu (U) > \mu (U +V_i - V_j ) > \mu (S + V_i ) $ 
\hfill  \framebox(7,7)

\section{Results on stability  and simplicity}

\begin{defin}
We say that a $G$-homogeneous  bundle is {\bf multistable}
 if it is the tensor product of a stable $G$-homogeneous bundle 
and an irreducible $G$-representation. 
\end{defin}

\begin{theorem} \label{RF} {\bf (Rohmfeld, Faini)} 
i)  \cite{Rohm}  A homogeneous bundle
$E$ is semistable if and only if 
 $\mu (F) \leq \mu(E)$ for any subbundle $F$ of $E$  induced by a 
subrepresentation  of the $P$-representation  inducing  $E$.

ii) \cite{Fa}  A homogeneous bundle
$E$ is multistable if and only if 
 $\mu (F) < \mu(E)$ for any subbundle $F$ of $E$  induced by a 
subrepresentation  of the $P$-representation  inducing  $E$.

\end{theorem}

We introduce now particular ${\cal Q}$-representations, called
 ``staircases''. Their importance is due to the fact that they are the 
${\cal Q}$-supports of  the homogeneous subbundles 
 of the homogeneous bundles whose
${\cal Q}$-supports are parallelepiped.

\begin{rema} \label{subquiv} { 
Let $E$ be a homogeneous bundle on ${\bf P}^n$ and $F$ be a homogeneous
 subbundle. 
Let   $S$ and $S'$  be the ${\cal Q}$-supports
of $E$ and   $F$ respectively.  By Theorem \ref{BKH} 
the ${\cal Q}$-representation of $F$ injects into 
the ${\cal Q}$-representation of $E$.   
If the multiplicities of $S$ are all $1$ 
and $S'$  contains the source of  an arrow $\beta$  in $S$ then 
$S'$ contains $\beta$. }
\end{rema}

\begin{defin} \label{scalini} We say that a 
 subgraph with multiplicities  of ${\cal Q} $ is  
a {\bf staircase} $S$  in a parallelepiped $R$ if all its multiplicities
 are $1$ and 
the  graph of $S$ is a subgraph of  $R$  satisfying the following property:
if $V$ is a vertex of $S$ then the  arrows of $R$ having $V$ as 
source must be arrows of $S$ (and then also their sinks must be vertices
 of $S$).

We say that a  subgraph  with multiplicities  
of ${\cal Q}$  is a staircase if it is a staircase in 
some parallelepiped.

\begin{center}
\includegraphics[scale=0.37]{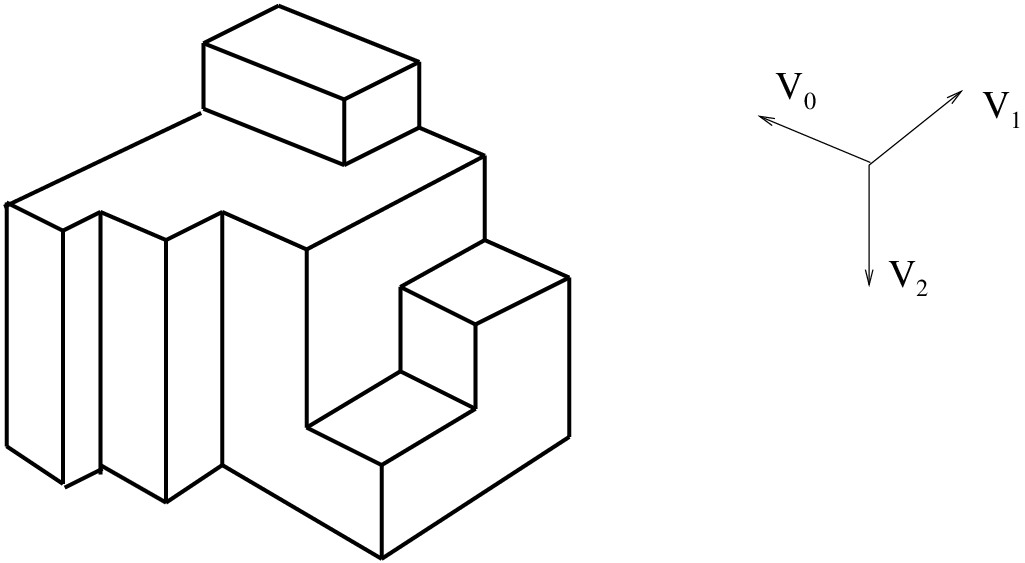}
\end{center}

Given a staircase $S$ in a  parallelepiped we define ${\cal V}_S$ to
 be the set of the vertices
of $S$ that are not sinks of any arrow of $S$.
We call the elements of ${\cal V}_S$ the {\bf vertices of the steps}. 
We say that a staircase has $k$ {\bf steps} if the cardinality of 
${\cal V}_S$ is $k$.

Let $V \in {\cal V}_S$.
We define the {\bf sticking out part} relative to  $V $ 
 as the part of $S$ whose vertices are exactly  the points  of $S$
greater than $V$ but not greater than any 
 other element of ${\cal V}_S$ and the arrows are 
all  the arrows connecting any of these vertices  (see Notation
\ref{multord}). 
\end{defin}

 By Remark \ref{subquiv} the ${\cal Q}$-support of a  homogeneous subbundle
 of a homogeneous bundle whose
${\cal Q}$-support is a  parallelepiped is a staircase in the parallelepiped.

\begin{rema} \label{prodtens} 
The support of $S^{l_1, l_2} Q(t) \otimes S^q V$ is the one shown in the 
figure, possibly cut off by  planes;
 precisely, if $q > l_1- l_2 $,  it is cut off  by a plane with direction
$\langle V_0, V_1 \rangle $ passing through the point 
$ \sigma \cap \{S^{l_1, l_2 }Q (t-q)  +k V_2 |\; k \in {\bf R} \} $ 
 and,  if $q >  l_2 $,  it is cut off  by  a plane with direction
$\langle V_1, V_2 \rangle $ passing through the point 
$ \pi \cap \{S^{l_1, l_2 }Q (t-q)  +k V_0 |\; k \in {\bf R} \} $, where
$\pi $ and $\sigma $ are the planes defined in  \S2.

(In fact, as $R$-representation,   $S^{l_1, l_2} Q(t) \otimes S^q V$
is equal to $S^{l_1, l_2} Q(t) \otimes ( \oplus_{i=0,..., q }  
 S^{q-i} Q (-i))$, apply Pieri's formula.)

Analogously $S^{l_1, l_2} Q(t) \otimes S^{q,q} V$ has the support
 shown in the figure below.
\end{rema}

\begin{center}
\includegraphics[scale=0.47]{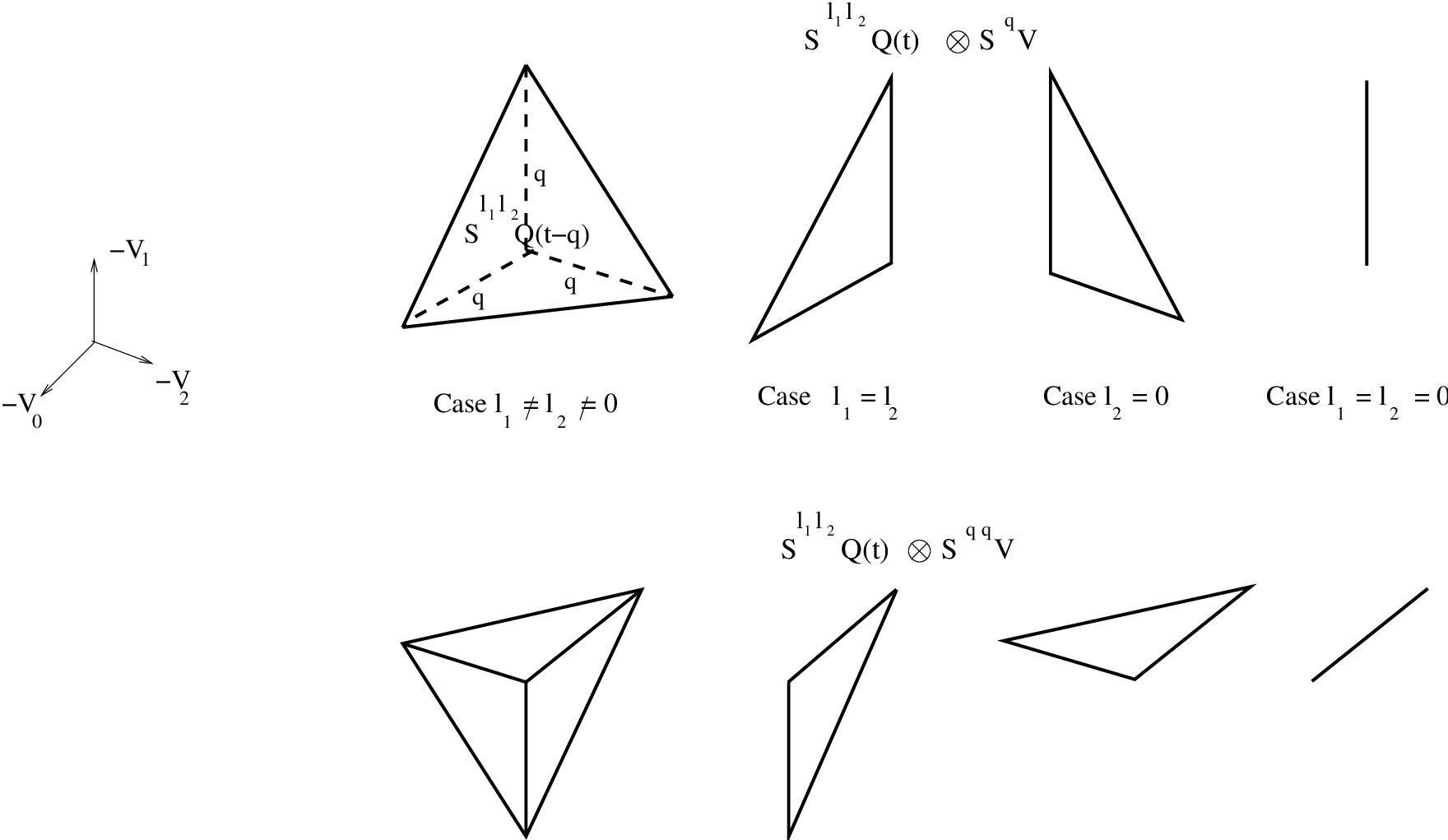}
\end{center}

\begin{theorem} \label{parstable} 
Let $E$ be the homogeneous vector bundle on ${\bf P}^3 = {\bf P}(V)$ 
whose ${\cal Q}$-support is a parallelepiped.  Then $E$ is stable 
(in particular it is simple) if it is not the tensor product of an 
$SL(V)$-representation and ${\cal O}(t) $ for some $t$ (that is, 
as we 
will see in \S4, a parallelepiped touching $\pi$ and $\sigma$).
\end{theorem}

{\it Proof.} First 
we prove that, to show that $E$ is stable, it is sufficient to
 show that it is multistable.

 If $E$ is  the tensor product of a stable homogeneous 
vector bundle $E'$  with an $SL(V)$-representation $W \neq {\bf C}$, then 
we can suppose that $W$ is irreducible; let $W = S^{p,q,r} V$.

$\bullet$ 
Suppose first that 
 $S^{p,q,r} V$, as $R$-representations, doesn't contain any $ {\cal O}(t)$;
then, as $R$-representations, the tensor product $E' \otimes S^{p,q,r} V$
is given by the tensor product of every summand of  $E'$ equal to 
$ {\cal O}(t)$ for some $t$ with $ S^{p,q,r} V$ (which is a parallelepiped)
and by the tensor products of every summand of  $E'$ different from
$ {\cal O}(t)$  $\forall t$ with every summand of  $ S^{p,q,r} V$ 
(which is a figure with more than one point and 
 parallel to $\langle V_1 -V_0 , V_2 -V_0 \rangle$).
 A union of such figures can't be a parallelepiped unless $E' =
 {\cal O}(t)$ for some $t$.

$\bullet$ 
Suppose now that $S^{p,q,r} V$, as $R$-representations, contains 
$ {\cal O}(t)$ for some $t$; then $q=r$. 

Observe that if $ p \neq q \neq 0$
 then   $S^{p,q,r} V$ would contain $S^{p-q-1} Q (-1) $ and  $S^{p-q+1,1} Q(-1)$
and so, if $  S^{l,m} Q(t)$ is not trivial, then $  S^{l,m} Q(t) \otimes
S^{p,q,r} V$ would contain points with multiplicity $2$.
Thus either $q=0$ or $p=q$. 

- If $ q=0 $ then, as $R$-representations, the tensor product $E' \otimes S^{p} V$
is given by the tensor product of all the  summands of  $E'$  equal to 
$ {\cal O}(t)$ for some $t$ with $ S^{p} V$ (which is a parallelogram
with sides parallel to $V_1 $ and $V_0 + V_1  + V_2$)
and by the tensor products of every summand of  $E'$ different from
$ {\cal O}(t)$  $\forall t$ with $ S^{p} V$ (figures shown in Remark 
\ref{prodtens}) and 
 a union of such figures can't be a parallelepiped unless $p=q=r=0$.

- Analogously the dual case  $p=q$.

\smallskip

To show that $E$ is multistable we consider the ${\cal Q}$-representation 
associated to $E$.

By Theorem \ref{RF}, $E$ is multistable if
 $\mu (F) < \mu(E)$ for any subbundle $F$ of $E$  induced by a 
subrepresentation  of the $P$-representation  inducing  $E$.
  Observe that, by Remark \ref{subquiv},
 the support of the ${\cal Q}$-representation  
of any such subbundle $F$ must  be a staircase  $C$ in 
 $R$ 
and vice versa any ${\cal Q}$-representation  whose support is a 
 staircase $C$ in $R$
 is the ${\cal Q}$-representation of a
subbundle  $F$  of $E$ induced by a subrepresentation  of the 
$P$-representation inducing  $E$.

We will show by induction on the cardinality  of ${\cal V}_C$ 
that $\mu(C) < \mu(R)$ for any $C$ staircase in $R$.

\underline{$k=1$ }
In this case $C$  is a subparallelepiped in the   parallelepiped $R$.
 Thus this case follows from Corollary \ref{rettangoli}.

\smallskip

\underline{$k-1 \Rightarrow k$} 
We will show that, given a staircase $C$  in $R$  with $k$ 
steps, there 
exists a  staircase  $C'$ in $R$ with $k-1$ steps such that  $\mu (C) \leq 
\mu(C') $. If we prove this, we conclude because $\mu(C)  \leq \mu (C') <
\mu(R) $, where the last inequality holds by  induction  hypothesis.

Let $C_1$ and $C_2$ be two staircases
 with $k-1$ steps 
obtained  from $C$ 
respectively ``removing and adding'' a parallelepiped $O$ and
a union $T$ of two parallepipeds.
($O$ is  a ``sticking out part'' of $C$ 
 and    $T$ is a nonempty union of parallelepipeds
 adjacent to  $O$ disjoint 
from $C$ such 
that the union of the points of $T$ with the point of $C$ gives
 a staircase with $k-1$ steps).

If $\mu (C_1)  \geq \mu (C) $ we conclude at once.

Thus we can suppose that $\mu(C_1 ) < \mu(C)$. We state that in this case 
$\mu (C_2 ) \geq  \mu (C)$.   In fact: let $\mu (C_1)= \frac{a}{b} $, 
$\mu(O) =\frac{c}{d}$ and $ \mu(T) =\frac{e}{f}, $ where the numerators 
are the first Chern classes and the denominators the ranks;  
since $\mu (C_1)  < \mu (C) $, we have $ \frac{a}{b} < \frac{a+c}{b+d}$, 
thus $ \frac{a}{b}  < \frac{c}{d}$; besides by Corollary \ref{rettangoli}
 $\mu(O) < \mu(T) $, 
i.e. $\frac{c}{d} < \frac{e}{f}$; thus 
$\frac{a+c+e}{b+d+ f} \geq  \frac{a+c}{b+d}$ i.e. $\mu (C_2 ) \geq  \mu (C)$.  
\hfill  \framebox(7,7)

\medskip

\begin{defin} We say that a  staircase contained in a  plane 
paralell to $\langle V_i, V_j \rangle$, for some $i,j$,  is 
{\bf completely
regular} if all the bundles corresponding the vertices of the steps 
(see Notation \ref{scalini}) are $P +l (V_i -V_j) $ $l=0,...,r$ for 
some $ P $ point of the quiver and $ r \in {\bf N}$.

We say that a staircase is a {\bf classical staircase} 
if it  is a cylinder on a completely regular staircase
in a  plane parallel to  $\langle V_1, V_2 \rangle $.
\end{defin}

\begin{theorem} \label{stairstable} 
Every classical staircase   is multistable
and it is stable unless it is either 
a cylinder of height $0$ (i.e. a completely regular staircase
in a  plane parallel to  $\langle V_1, V_2 \rangle $) or
 a  parallelepiped given by the tensor product of an 
$SL(V)$-representation and ${\cal O}(t) $ for some $t$.   
\end{theorem}

\begin{center}
\includegraphics[scale=0.39]{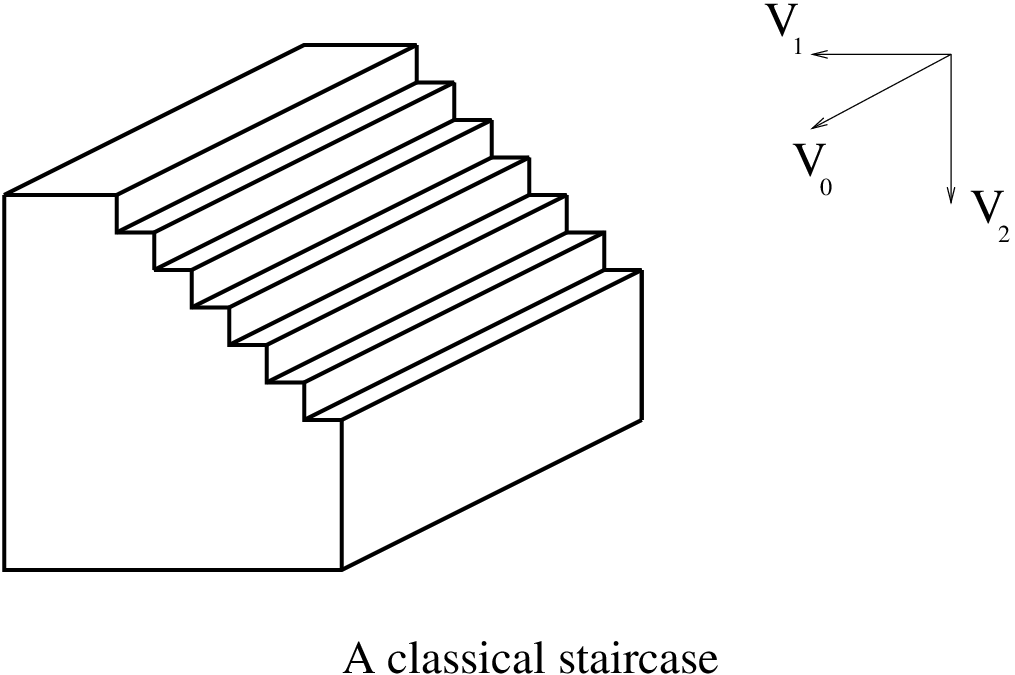} 
\end{center}

{\it Proof.} First we prove the statement on multistability.

{\sf Fact 1.} {\it  For any classical staircase, let's order
the vertices $P_1,P_2,...$ by going in the direction $ V_1 -V_2$
and let $R_i$ be the parallelepiped contained in the staircase 
whose vertices are exactly  those greater than $P_i$. Let $ H_i =
R_i -R_{i-1}$, $ R_0 = \emptyset $ (horizontal steps)
and  $ E_i =
R_i -R_{i+1}$, $ R_0 = \emptyset $  (vertical steps).
 We have 
$$\mu (H_i)  > \mu (H_{i-1}) \hspace{2cm}  \mu (E_i ) > \mu (E_{i+1})$$
for any $i$.}

\smallskip

{\sf Proof.}
It follows from Proposition \ref{main} and Corollary \ref{rettangoli}. 

\smallskip

{\sf Fact 2.} {\it 
Let $S$ be a classical staircase.
 Then for every sticking out part $O$ of $S$
we have $$ \mu (O) > \mu (S-O)$$
Therefore  $$\mu (S) > \mu (S-O) $$
More generally let us define a ``piece $O$ of the staircase $S$'' 
in the following way. Let $P$ and $Q$ be two vertices with $Q =  P
+ m (V_2 -V_1)$ and let us consider the triangle $R$ with vertices 
$P,Q, Q +m V_1$  
and let $O$ be the part of the staircase cylinder on $R$ (see the figure 
below part a, 
section for a plane parallel to $ \langle V_1 ,V_2  \rangle $ ).

If $  O_1,..., O_k$ are pieces of the staircases $S$, we have that 
$\mu (O_i) > \mu (S -O_1 ...-O_{i-1} -O_{i+1} ...-O_k )$. }

\smallskip

{\sf Proof.}
Let $b$ be the  plane on which the base of $O$ is and  
let $l$ be the  plane on which the left side of $O$ is. 
Let $T_1$ be the 
staircase  whose vertices are the vertices of $S$ that are  
either above $b$ or on $b$ and  on the left  of $l$ (see the figure 
below part b, section for a plane parallel to $ \langle V_1 ,V_2  \rangle $ ).
Let $T_2$ be the 
staircase  whose vertices are the vertices of $S$ that  are  
below  $b$ and either on the right of $l$  or on   $l$.

Let $K$ be the rectangle $$K= S -T_1 -T_2 -O$$

\begin{center}
\includegraphics[scale=0.36]{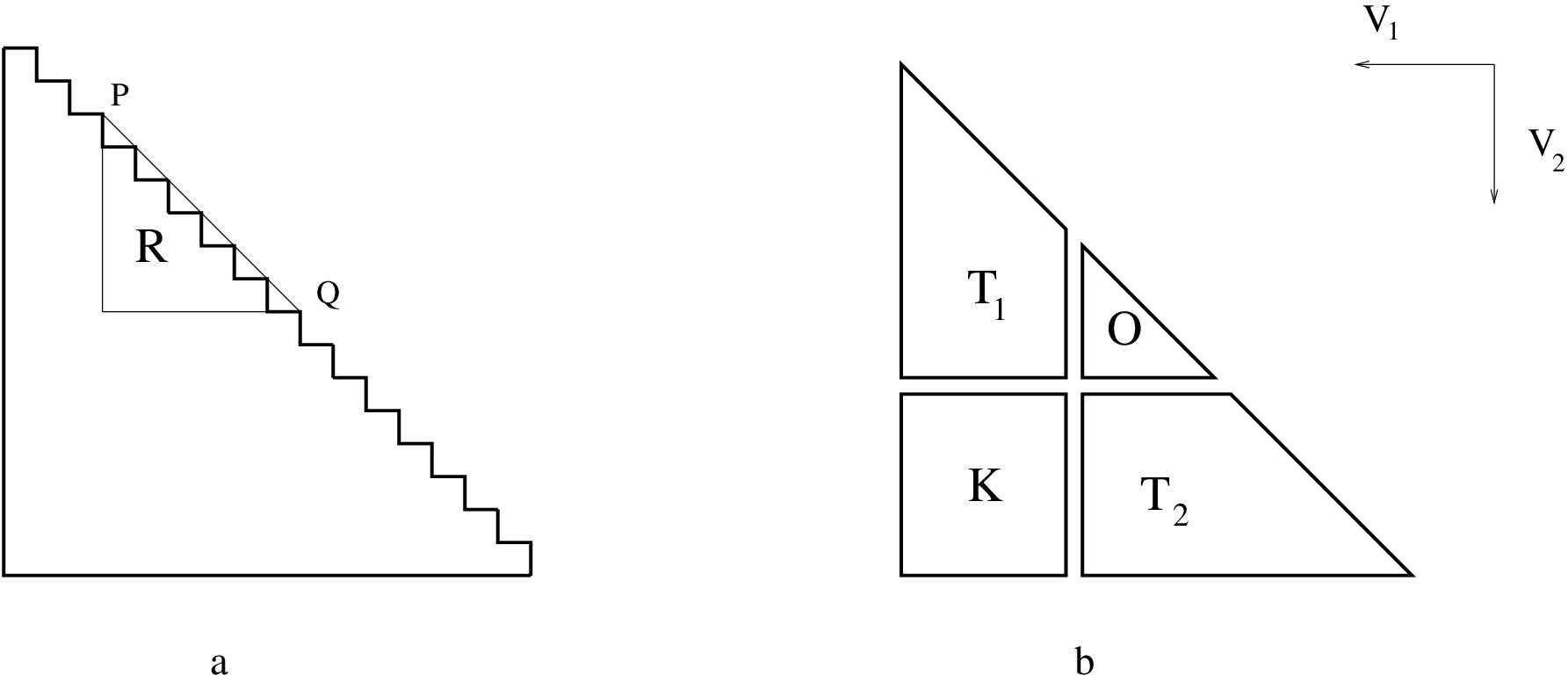}
\end{center}

By Corollary \ref{rettangoli} $\mu (O) > \mu (K) $.
Besides, by applying Fact 1 to the staircases $T_1 +O$ and $T_2 +O$
(where $T_i+O$ is the smallest staircase containing $T_i$ and $O$), we get 
 $$\mu (O) > \mu (T_1) \;\;\;\;\;\;\;\;\;\; \mu (O) > \mu (T_2)$$ 
Hence $\mu(O) > max \{\mu (K) , \mu(T_1) , \mu(T_2)\} \geq \mu(S-O)$ 
(see Remark \ref{gio}). So we have proved the first statement of Fact 1.

Analogously for the second statement of Fact 2.

\smallskip

Now we are ready to prove that every bundle  such that  its  ${\cal Q}$-support
is a classical 
staircase  $S$   is multistable.
Let $C$ be the support 
of a ${\cal Q}$-representation subrepresentation of $S$ (thus again a 
staircase by Remark \ref{subquiv}). We want to prove $\mu(C) < \mu (S) $  by 
induction on the cardinality $k$ of $\nu_C$.

\underline{$k=1$}. The statement follows from
Corollary \ref{rettangoli} and Fact  1.

\underline{$k-1 \Rightarrow k$.} To prove this implication 
 we do induction on $$- area
(bd(C) \cap bd (S)) $$ 
where $bd $ denotes the border and the border of a staircase
 is the border of the part of the space inside the staircase. 

Let $C$ be a staircase with $k$ steps  support 
of a subrepresentation of $S$.

$\bullet$
If $\mu (C -O) \geq \mu (C) $ for some sticking out part $O$ of $C$, we 
conclude at   once because $C-O $ has $k-1$ steps; thus by induction 
assumption $\mu(S) > \mu (C-O)$ and then $\mu (S) > \mu (C)$.  

$\bullet$ 
Thus we can suppose  $\mu (C) > \mu (C-O) $ for every sticking out
 part $O$ of $C$ 
i.e.  $\mu (O) > \mu (C) $ for every sticking out part $O$ 
of $C$.

$- $
Suppose there exists a sticking out part $O$ of $C$ 
such that there exists $A$ parallelepiped or union of two
parallepipeds such that  $A$ is disjoint from $C$,
 $A+C$ is a staircase with less steps
(where $A+C$ is the smallest staircase containing $C$ and $A$), a side
 of $A$ is equal to a side of $O$  and
$A':=A \cap S \neq \emptyset$.

Since $\mu (A') > \mu (O) $ by Corollary \ref{rettangoli}
 and $\mu(O)>  \mu (C) $ by assumption,
 we have $\mu (A') > \mu (C)$ and thus
  \begin{equation} \label{st} \mu (C+ A') > \mu (C) \end{equation} 
If $A$ is a subgraph of $S$ i.e. $A= A'$, then $C + A' $ is a 
staircase  with $k-1$ steps thus, by induction assumption, $\mu (S) > \mu
 (C+ A') $; hence 
 $\mu (S) > \mu (C)$  by (\ref{st}).

If $A$ is not a subgraph of $S$ i.e. $A\neq  A'$, then $area
 (bd (C+A') \cap bd (S)  ) > area (bd(C) \cap bd (S))$ and by induction
 assumption $\mu (S) > \mu (C +A') $; hence we conclude again
$\mu (S) > \mu (C)$  by (\ref{st}).

$-$ Suppose there doesn't exist a sticking out part $O$ of $C$ 
such that there exists $A$ parallelepiped or union of two
parallepipeds such that  $A$ is disjoint from $C$,
 $A+C$ is a staircase with less steps, a side of 
$A$ is equal to a side of $O$  and
$A':=A \cap S \neq \emptyset$,
that is for every  sticking out part $O$ of $C$ and for every
$A$ parallelepiped or union of two
parallepipeds such that  $A$ is disjoint from $C$,
 $A+C$ is a staircase with less steps and  a side of 
$A$ is equal to a side of $O$, we have 
$A':=A \cap S =\emptyset$.
Observe that in this case $C$ must be a cylinder on a staircase 
in a plane parallel to $ \langle V_1, V_2\rangle$. 
Then there exists a chain of staircases 
 $C = C_0 \subset C_1 \subset ... \subset C_r =C'$
 such that  $C_i $ is obtained 
from $C_{i+1}$ taking off one of its sticking out parts and $C'$ is 
the intersection of $S$ with  a semispace whose border plane 
is parallel 
to $\langle V_1, V_2\rangle$; thus, by Fact 2, $\mu (C) \leq \mu (C')$ 
and we can prove that $\mu (C') \leq \mu (S)$  in an analogous way 
as Proposition  \ref{main}; so we 
conclude the proof of the statement on multistability.

As to stability, the proof is completely analogous to the proof 
in Theorem \ref{parstable} \hfill  \framebox(7,7)

\section{Resolutions of parallelepipeds and staircases}

In this section we investigate the minimal free resolutions of the bundles
 whose quiver supports are parallelepipeds or staircases.

Firstly observe that the supports of the bundles 
$S^{\lambda_1, \lambda_2, \lambda_3} V (t)$
are parallelepipeds with an edge on $\pi$ and an edge on $\sigma$ 
(border planes of the quiver, they are defined in \S4).

\begin{rema}
 Let $R$ be a parallelepiped.  
We can get the minimal free resolution of the bundle relative to $R$ 
in the following way: 
let  $Q,S,P,T$ be the parallepipeds as shown in the figure
(that is: let $Q$ be the parallelepiped touching $\sigma$ 
and whose edges in the directions $V_0$ and $V_2$ 
have the same length of the corresponding edges of $P$ 
(we get it by going in the direction  $V_1$) and so on). We denote
by $R \cup Q$ the minimum parallelepiped containing both $R$ and $Q$ and
analogously for the others.

 We get the minimal free resolution (where we identify the bundles with their
${\cal Q}$-supports)
$$ 0 \rightarrow T \cup S \rightarrow (T \cup S \cup P) 
 \oplus (Q \cup S)  \rightarrow R \cup P \cup Q \cup S \rightarrow R
 \rightarrow 0 $$  (all the components of the maps 
nonzero).
Thus the resolution of a bundle $R$ whose support is 
 a parallelepiped   touching neither $\pi$ nor $ \sigma $ is
$$ 0 
\rightarrow 
S^{\lambda_1-l, \lambda_2-k, \lambda_3 } V (t-k-l) 
\rightarrow 
S^{\lambda_1, \lambda_2-k, \lambda_3 } V (t-k) 
\oplus 
S^{\lambda_1-l, \lambda_2, \lambda_3 } V (t-l) \rightarrow
$$ $$\rightarrow 
S^{\lambda_1, \lambda_2 , \lambda_3 }
 V (t)
\rightarrow 
R 
\rightarrow 
0 $$
(with $\lambda_1, \lambda_2, \lambda_3 , t,k,l \in {\bf N}$,  $ k,l
 \geq 1$  and 
$\lambda_1 \geq \lambda_2 \geq \lambda_3 $ and all the components
 of the maps nonzero).

If $R$ touches $\pi$, then $T$, $S$  and $P$ 
are missing and the minimal free 
resolution becomes 
$$ 0 
\rightarrow 
S^{\lambda_1-l, \lambda_2, \lambda_3 } V (t-l) \rightarrow
S^{\lambda_1, \lambda_2 , \lambda_3 }
 V (t)
\rightarrow 
R 
\rightarrow 
0 $$
If $R$ touches $\sigma$, then $Q$ and $S$ are missing and the minimal free 
resolution becomes 
$$ 0 
\rightarrow 
S^{\lambda_2- 1, \lambda_2-k, \lambda_3 } V (t-k- \lambda_1 +\lambda_2-1) 
\rightarrow 
S^{\lambda_1, \lambda_2-k, \lambda_3 } V (t-k) 
\rightarrow 
S^{\lambda_1, \lambda_2 , \lambda_3 }
 V (t)
\rightarrow 
R 
\rightarrow 
0 $$

\end{rema}

\begin{center}
\includegraphics[scale=0.39]{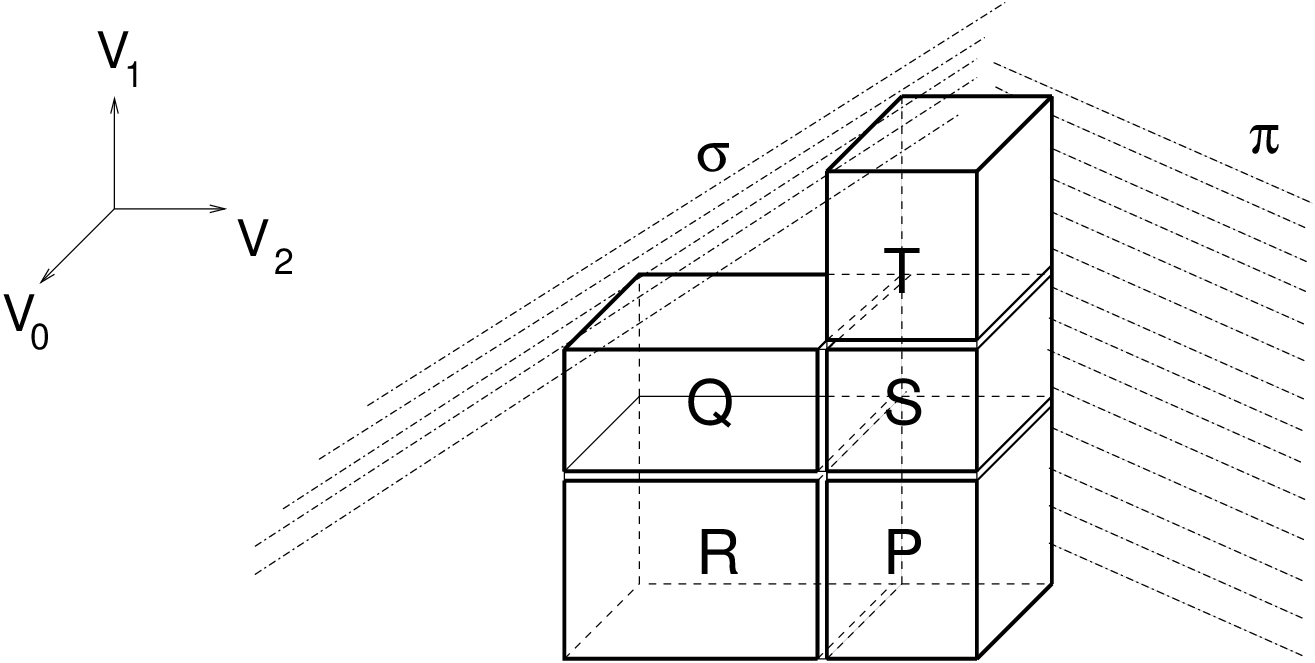}
\end{center}

\medskip

The  above remark shows that 
Theorem \ref{parstable} can be restated in  Theorem \ref{parstablebis}.

\medskip

Now we want  to study the resolution of the bundle corresponding to
a staircase which is a cylinder on a staircase
in a  plane parallel to  $\langle V_1, V_2 \rangle $. 
For any such staircase $S$ let's order
the vertices $P_1,P_2,...$ by going in the direction $ V_1 -V_2$
and let $R_i$ be the parallelepiped contained in the staircase 
whose vertices are exactly  those greater than $P_i$. Let  $ E_i =
R_i -R_{i+1}$, $ R_0 = \emptyset $  (vertical steps).

First let's suppose that $S$ touches $\pi$.
Let $K_i$ be the parallelepiped  touching $\sigma$ and $\pi$ 
containing $E_i$ and with the length of the  edges 
in the direction $V_0 $ and $V_2$ 
equal to the corresponding edges of $E_i$ 
(we get it by ``going''in the 
direction of $V_1$).

We get the minimal free resolution:
$$ 0 \rightarrow \oplus_i (K_i -E_i) 
\rightarrow \oplus_i K_i \rightarrow S \rightarrow 0$$
where the second map restricted 
to $K_i- E_i$ has only the components  $K_i -E_i \rightarrow K_i$ 
and  $K_i -E_i \rightarrow K_{i+1}$ nonzero.  

Now let's suppose that $S$ doesn't touch $\pi$.
Let $E'_i$ be the parallelepiped containing $E_i$ 
and with the length of the  edges 
in the direction $V_0 $ and $V_1$ 
equal to the corresponding edges of
 $E_i$ (we get it by ``going''in the 
direction of $V_2$). 
Let $K'_i$ be the parallelepiped  touching $\sigma$ and $\pi$ 
containing $E'_i$ and with the length of the edges 
in the direction $V_0 $ and $V_2$ 
equal to the corresponding edges of $E'_i$ (we get it by ``going''in the 
direction of $V_1$). 

Let $Z$ be the parallelepiped  touching $\sigma$ and $\pi$ 
containing $E'_1 -E_1$ and with the length of the edges 
in the direction $V_0 $ and $V_2$ 
equal to the corresponding edges of $E'_1 - E_1$ 
(we get it by ``going''in the 
direction of $V_1$).

We get the minimal free resolution:
$$ 0 \rightarrow Z - \oplus_i  E'_i \rightarrow Z \oplus 
\oplus_i (K'_i -E'_i) 
\rightarrow \oplus_i K'_i \rightarrow S \rightarrow 0$$
where the third map restricted 
to $K'_i- E'_i$ has only the components  $K'_i -E'_i \rightarrow K'_i$ 
and  $K'_i -E'_i \rightarrow K'_{i+1}$ nonzero.

\begin{center}
\includegraphics[scale=0.37]{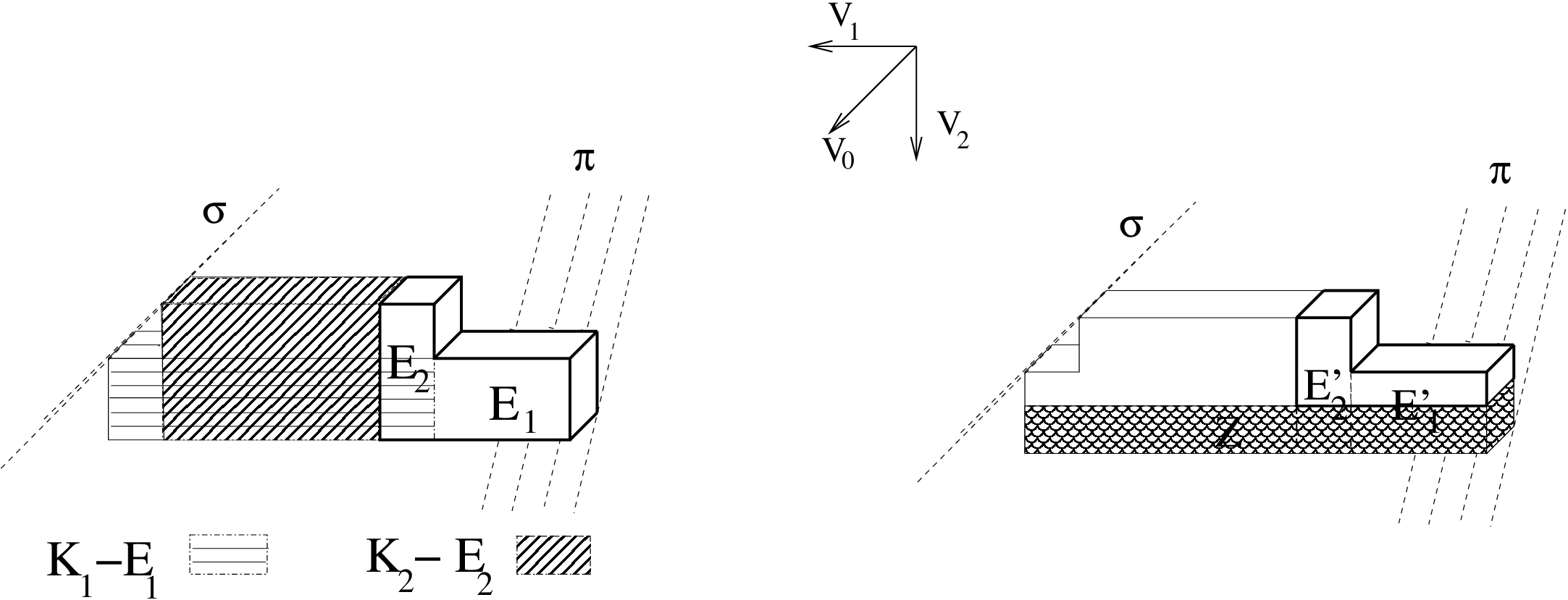}
\end{center}

So we have that the resolution of the bundle $E$ corresponding to
a staircase which is a cylinder on a  staircase
in a  plane parallel to  $\langle V_1, V_2 \rangle $ and whose vertices
are lined up in a line parallel to $ V_1-V_2$ 
 is, 
for some $\lambda_1, \lambda_2, \lambda_3 , r, s_1,...,s_r, t_1,...,t_r
 \in {\bf N}$, $ s_1>...>s_r$, $s_i = s_{i+1} + t_{i+1}$ for $i=1,....,
r-1$, $\epsilon,
 \delta
\in \{0,1\}$
$\lambda_1 \geq \lambda_2 \geq \lambda_3 $ and  $SL(V)$-invariant maps:
{\footnotesize
$$ 0 
\rightarrow 
\delta S^{\lambda_1 +s_r -k, \lambda_2 -s_1 -k, \lambda_3 } V (-2k +s_r
-s_1) 
\stackrel{\varphi}{\rightarrow}
\delta S^{\lambda_1+s_1 -k, \lambda_2 -s_1-k, \lambda_3 } V (-2k) 
\oplus \oplus_{i=1,...,r-\epsilon}
S^{\lambda_1 +s_i, \lambda_2-s_i, \lambda_3 } V  \stackrel{\psi}{
\rightarrow}
$$ $$\rightarrow \oplus_{i=1,...,r}
S^{\lambda_1+s_i+t_i, \lambda_2 -s_i, \lambda_3 }
 V (t_i)
\rightarrow 
E
\rightarrow 
0 $$}
where the only nonzero component of $\varphi$  is the first,
  $\psi|_{S^{\lambda_1 +s_i, \lambda_2-s_i, \lambda_3 } V} $ 
has only the components  into 
$ S^{\lambda_1+s_i+t_i, \lambda_2 -s_i, \lambda_3 } V (t_i) $ and   
$S^{\lambda_1+s_{i+1}+t_{i+1}, \lambda_2 -s_{i+1}, \lambda_3 } V (t_{i+1})$
 nonzero and  
$\psi|_{S^{\lambda_1+s_1 -k, \lambda_2 -s_1-k, \lambda_3 } V (-2k)} $
 has all the components nonzero.
($\delta$ is $0$ iff the staircase touches $ \pi $ and 
$\epsilon$ is $1$ iff the staircase touches $ \sigma $).

\bigskip

So Theorem \ref{stairstable} can be restated in
 Theorem \ref{stairstablebis}.

\section{ Proof of  Theorem \ref{kerscalasempl}}

\begin{lemma} \label{lemmaL}
Let $\lambda_1,..., \lambda_n,s \in {\bf N}$  with $\lambda_1
\geq ...\geq  \lambda_n$.

Let $T=  \{(s_1, ....,s_{n+1}) \in {\bf N}^{n+1}
 |\; s_1 +.....+s_{n+1} =s, \; s_i \leq 
\lambda_i- \lambda_{i+1} \;\; i=2,..., n
\;\;  s_{n+1} \leq \lambda_n \}$. For every $ M \subset T$
let ${\cal P}_M$ be the following statement: 
 for every  $V$  complex vector space of dimension $n+1$, the 
commutativity of the diagram  of bundles
on ${\bf P}(V)$ {\small $$\begin{array}{ccc}  S^{ \lambda_1,..., \lambda_n
 } V (-s) &
\stackrel{\varphi}{\longrightarrow} & \oplus_{(s_1,....,s_{n+1}) \in M}
 S^{ \lambda_1+s_1,....,\lambda_n +s_n, s_{n+1} } V   \\ {\scriptsize A}
\downarrow \;\;\;\;\;\;\; & & \;\;\; \downarrow {\scriptsize B}
\\ S^{ \lambda_1,..., \lambda_n
 }V (-s) & \stackrel{\varphi}{\longrightarrow} 
 & \oplus_{(s_1,...,s_{n+1}) \in M} 
S^{\lambda_1+s_1,....,\lambda_n +s_n, s_{n+1}} V
  \end{array}$$ } \hspace*{-0.3cm} 
(where $A$ and $B$ are linear maps and the components of $\varphi$  are
nonzero $SL(V)$-invariant maps) implies  
$A= \lambda I$ and $B = \lambda I$ for some $\lambda \in {\bf C}$.

Let $C \subset T$ with $C \neq \emptyset $, $T-C \neq \emptyset$. 
Then ${\cal P}_C$ is true if and only if ${\cal P}_{T-C}$ is true.
\end{lemma}

{\it Proof.}  Completely analogous to the proof of Lemma 38 in \cite{O-R1}.

\bigskip

{\it Proof of Theorem \ref{kerscalasempl}.} Let 
$T= \{s \in {\bf N}^4 |
\; \; 
s_i \leq \lambda_{i-1} -\lambda_i \; for\; i=2,3 ,\;\;
s_4 \leq \lambda_{3},\;\;  s_1+s_2+s_3+s_4=c\}$ and $C = 
\{s \in T |
s_1=d \;\; s_4=0\}$.
Thus $K =T-C$. 

 First let us suppose that $ \lambda_3 \neq 0$. 
Let $f$ be an endomorphism of $E$.  
It induces a commutative diagram
$$ \begin{array}{ccccccc} 0 \rightarrow 
& S^{\lambda_1,\lambda_2,\lambda_3 } V
& \stackrel{\alpha}{\rightarrow} & \oplus_{s \in T- C} 
S^{\lambda_1+s_1,\lambda_2+s_2,\lambda_3+s_3 ,s_4} V (s_1+s_2+s_3+s_4) &
\rightarrow  & E & \rightarrow 0
\\ & A \downarrow & & \downarrow B & &\downarrow f  &
\\ 0 \rightarrow & S^{\lambda_1,\lambda_2,\lambda_3 } V
& \stackrel{\alpha}{\rightarrow} & \oplus_{s \in T- C} 
S^{\lambda_1+s_1,\lambda_2+s_2,\lambda_3+s_3, s_4 } V (s_1+s_2+s_3+s_4) &
\rightarrow  & E & \rightarrow 0
\end{array} $$
(In fact:
 write $ 0 \rightarrow R \rightarrow S \rightarrow E \rightarrow 0$ 
for the minimal resolution of $E$ for short;
 applying $Hom (S, \cdot ) $ to it,
 we can prove the existence and uniqueness of $B$ 
and by applying $Hom (R, \cdot ) $ we can prove the existence of $A$).
We have to prove that $B=I$.

By Lemma \ref{lemmaL}, to show our statement, it is sufficient 
to show that the vertical map 
of
a diagram \begin{equation} \label{f1} \begin{array}{ccc} 
 S^{\lambda_1,\lambda_2,\lambda_3 } V
& \stackrel{\varphi}{\rightarrow} & \oplus_{s \in C} 
S^{\lambda_1+s_1,\lambda_2+s_2,\lambda_3+s_3 } V (s_1+s_2+s_3) 
\\ \downarrow D & & \downarrow B' 
\\  S^{\lambda_1,\lambda_2,\lambda_3 } V
& \stackrel{\varphi}{\rightarrow} & \oplus_{s \in C} 
S^{\lambda_1+s_1,\lambda_2+s_2,\lambda_3+s_3 } V (s_1+s_2+s_3) 
\end{array} \end{equation}
(where the horizontal arrows are $SL(V)$-invariant) are the identity maps.

Observe that   $D(Ker (\varphi)) \subset  Ker (\varphi)$.
Thus we have a commutative diagram 
 \begin{equation} \label{f2} \begin{array}{ccc}  Ker (\varphi) &
\longrightarrow &  S^{\lambda_1,\lambda_2,\lambda_3 } V
 \\ \;\;\;\;\;\;\;\;\;\;\;\downarrow
D|_{Ker(\varphi)}  & &  \downarrow D 
\\  Ker (\varphi) & \longrightarrow & S^{\lambda_1,\lambda_2,\lambda_3 } V
\end{array}
\end{equation}
Let $ 0 \rightarrow P \rightarrow Q \rightarrow U
 \rightarrow Ker (\varphi) \rightarrow 
0 $
 be a minimal free resolution of $Ker(\varphi)$. 
Since $Ker(\varphi)$ is simple (it is a classical staircase  with  nonzero
 height since $\lambda_3 \neq 0$) 
the map $D|_{Ker(\varphi)} : Ker (\varphi) \rightarrow Ker (\varphi)$ 
is the identity.
Thus we get a commutative diagram \begin{equation} \label{f3}
 \begin{array}{ccc} Q 
 & \rightarrow  & Ker(\varphi) 
\\  \downarrow  I & & \;\;\;\;\downarrow  D|_{Ker(\varphi)}= I \\ 
Q & \rightarrow &  Ker (\varphi) 
\end{array}
\end{equation}
and then
\begin{equation} \label{f4}
\begin{array}{ccc} Q_{max}  & \rightarrow  & Ker(\varphi) 
\\  \downarrow  I & & \;\;\;\;\downarrow D|_{Ker(\varphi)}=I  \\ 
Q_{max}  & \rightarrow &  Ker (\varphi) 
\end{array} \end{equation}
where $Q_{max}$ is the sum of the addenda of $Q$ with maximal twists,
and then (composing (\ref{f2}) with (\ref{f4})) 
 \begin{equation} \label{f5}
 \begin{array}{ccc}
 Q_{max}  & \rightarrow  &  S^{\lambda_1,\lambda_2,\lambda_3 } V  
\\  \downarrow  I & & \;\;\;\;\downarrow D  \\ 
Q_{max}  & \rightarrow &   S^{\lambda_1,\lambda_2,\lambda_3 } V
\end{array} \end{equation}
Thus $D=I$ (consider $H^0( \cdot^{\vee})^{\vee}$ of all the maps of 
(\ref{f5})).

If $ \lambda_3 =0 $ the proof is completely analogous to the one 
of Lemma 45 in \cite{O-R1}.
\hfill  \framebox(7,7)


\end{document}